\documentclass[a4paper,10pt]{amsart}

\usepackage{amsfonts,amsthm,amssymb,amsmath,amscd,ascmac}

\usepackage{latexsym}
\usepackage{mathrsfs}

\def\qed{\hfill $\Box$}
 
\newcommand{\DS}{\displaystyle}

\newcommand{\NN}{\mathbb{N}}
\newcommand{\ZZ}{\mathbb{Z}}
\newcommand{\QQ}{\mathbb{Q}}
\newcommand{\RR}{\mathbb{R}}
\newcommand{\CC}{\mathbb{C}}
\newcommand{\PP}{\mathbb{P}}
\newcommand{\FF}{\mathbb{F}}

\newcommand{\ee}{\mathbf{e}}
\newcommand{\z}{\zeta}
\newcommand{\fA}{\mathfrak{A}}
\newcommand{\dv}{\mathrm{div}}

\newcommand{\HH}{\mathbb{H}}
\newcommand{\DD}{\mathbb{D}}

\newcommand{\I}{\mathrm{I}}

\begin{document}
\title{The Fermat septic and the Klein quartic as moduli spaces of hypergeometric 
Jacobians}
\date{}
\author{Kenji Koike}  
\dedicatory{Dedicated to the 70th birthday of Professor Hironori Shiga.}
\address{Faculty of Education \\
University of Yamanashi \\
Takeda 4-4-37, Kofu \\
Yamanashi 400-8510, Japan\\}
\date{}
\email{kkoike@yamanashi.ac.jp}
\begin{abstract}
We give uniformizations of the Klein quartic curve and the Fermat septic 
curve as Shimura curves parametrizing Abelian $6$-folds with endomorphisms 
$\ZZ[\z_7]$. 
\end{abstract}
\maketitle
\section{Introduction}
The Gauss hypergeometric differential equation 
\[
  E(a,b,c) \ : \ z(z-1) u'' + \{(a+b+1)z -c \} u' + ab u = 0
\]
is regular on $\CC-\{0,1\}$ for general parameters $a, b$ and $c$, and the solution 
space is spanned by Euler type integrals
\[
 \int_{\gamma} x^{a-c} (x-1)^{c-b-1} (x-z)^{-a} dx,
\]
that are regarded period integrals for algebraic curves if $a, b, c \in \QQ$. 
Two independent solutions $f_0(z), f_1(z)$ define a multi-valued analytic function
$\mathfrak{s}(z) = f_0(z) / f_1(z)$ (Schwarz map), and monodromy transformations for 
$\mathfrak{s}(z)$ are given by fractional linear transformations. 
\\ \indent
If parameters satisfy the conditions
\[
 |1-c| = \frac{1}{p}, \quad |c-a-b| = \frac{1}{q}, \quad |a-b| = \frac{1}{r},
\qquad \frac{1}{p} + \frac{1}{q} + \frac{1}{r} <1,
\]
with $p,q,r \in \NN \cup \{ \infty\}$, the monodromy group is isomorphic 
to a triangle group
\[
 \Delta(p,q,r) = \left< M_0, \ M_1, \ M_{\infty} \ | \ 
 M_0^p = M_1^q = M_{\infty}^r = M_0 M_1 M_{\infty} = 1 \right>
\]
(the condition $M_0^p = 1$ is omitted if $p = \infty$, and so on). In this case, 
the upper half plane is mapped to a triangle with vertices 
$\mathfrak{s}(0), \ \mathfrak{s}(1)$ and $\mathfrak{s}(\infty)$, angles 
$\pi /p, \ \pi/q$ and $\pi /r$ respectively by $\mathfrak{s}$, and so is the lower 
half plane. Copies of these two triangles give a tessellation of a disk $\DD$ by 
the monodromy action, and we have an isomorphism 
$\overline{\DD/\Delta(p,q,r)} \cong  \overline{\CC-\{0,1\}} = \PP^1$.
For example, $E(1/2,1/2,1)$ is known as the Picard-Fuchs equation for the Legendre family 
of elliptic curves $y^2 = x (x-1)(x-z)$ and the monodromy group 
$\Delta(\infty, \infty, \infty)$ is projectively isomorphic to the congruence subgroup 
$\Gamma(2)$ of level $2$. Also a triangle group $\Delta(n,n,n)$ with $n \geq 4$ is 
interesting, since 
its commutator subgroup $N_n$ gives a uniformization of the Fermat curve $\mathcal{F}_n$ 
of degree $n$. More precisely, the natural projection 
$\DD / N_n \rightarrow \DD/ \Delta(n,n,n) = \PP^1$ is an Abelian covering 
branched at $0, 1$ and $\infty$ with the covering group 
$\Delta(n,n,n) / N_n \cong (\ZZ /n \ZZ)^2$ (see \cite{CIW94}).     
\\ \indent
In \cite{T77}, Takeuchi determined all arithmetic triangle groups. According to it, 
$\Delta(n,n,n)$ is arithmetic (and hence the Fermat curve $\mathcal{F}_n$ is a Shimura 
curve) for $n \in FT = \{4, 5, 6, 7, 8, 9, 12, 15 \}$. These groups come from 
the Picard-Fuchs equation for algebraic curves $X_t : y^m = x(x-1)(x-t)$ with 
$m = n$ (resp. $m=2n$) if $n \in FT$ is odd (resp. even). 
Among them, $n=5$ and $7$ are special in the sense that a Jacobian $J(X_t)$ is simple in 
general, and Picard-Fuchs equations describe variations of Hodge structure on 
the whole of $H^1(X_t, \QQ)$, rather than sub Hodge structures. These two families are 
treated by Shimura as examples of PEL families in \cite{Sm64}.
Also de Jong and Noot studied them as counter examples of Coleman's conjecture 
(which asserts the finiteness of the number of CM Jacobians for a fixed genus $g \geq 4$) 
for $g=4, 6$ in \cite{dJN91} (see also \cite{R09} and \cite{MO13} for this direction).
\\ \indent
For $n=5$, we gave $\mathfrak{s}^{-1}$ by theta constans in \cite{K03} as a byproduct of 
study of the moduli space of ordred five points on $\PP^2$.
In present paper, we compute the monodromy group, Riemann's period matrices and the 
Riemann constant with an explicit symplectic basis for $n=7$. Using them, we express 
the Schwarz inverse map $\mathfrak{s}^{-1}$ by Riemann's theta constants 
(Theorem \ref{main}). As a consequence, we give explicit moduli interpretations of the Klein 
quartic curve $\mathcal{K}_4$ and the Fermat septic curve $\mathcal{F}_7$ as modular varieties 
parametrizing Abelian $6$-folds with endomorphisms $\ZZ[\z_7]$. 
(Corollary \ref{Klein1}, Corollary \ref{Klein2}). The Klein quartic is classically 
known to be isomorphic to the elliptic modular curve of level $7$. In \cite{E99}, Elkies 
studied it as a Shimura curve parametrizing a family of QM Abelian $6$-folds. Our 
interpretation of $\mathcal{K}_4$ gives the third face as a modular variety.  
Our expression of $\mathfrak{s}^{-1}$ is a variant of Thomae's formula. This kind of formula 
for cyclic coverings was studied in general context by Bershadsky-Radul 
(\cite{BR87}, \cite{BR88}), Nakayashiki (\cite{Na97}) and Enolski-Grava (\cite{EG06}), 
but our standpoint is more moduli theoretic as a classical work of Picard 
(\cite{P1883}) which produces modular forms on a $2$-dimensional complex ball. 
In \cite{Sh88}, Shiga determined Picard modular forms explicitly, and his results were 
applied to number theory and cryptgraphy (see \cite{KS07} and \cite{KW04}). We expect that 
also our concrete results will give a good example to develop a generalization of arithmetic 
theory of elliptic curves. Here we mention that there are several studies of automorphic 
forms for triangle groups (e.g. \cite{Mi75}, \cite{W81}, \cite{H05} and \cite{DGMS13}). 
However explicit constructions of autmorphic forms for co-compact triangle groups in the 
view point of the Picard's work seems to be very few. 
\\ \indent
Our Schwarz map is regarded also as a periods map of K3 surfaces. In pioneer work 
\cite{Sh79-81}, Shiga studied families of elliptic K3 surfces with period maps to 
complex balls. These K3 surfaces have a non-symplectic automorphism of 
order $3$, which induces a Hermitian structure on the transcendental lattice. Now K3 
surfaces with non-symplectic automorphisms of prime order are classified (see \cite{AST11}), 
and many of them are known to be quotients of product surfaces (\cite{GP}). 
In the last section, we give elliptice K3 surfaces $S_t$ associated to $X_t$ and compute the 
Neron-severi group and the Mordell-Weil lattice of $S_t$.  
\section{Uniformization of Fermat Curves}
\subsection{Hypergeometric integral} 
We compute monodromy groups and invariant Hermitian 
forms for hypergeometric integrals
\[
u(t) = \int \Omega_{\alpha}(x), \qquad \Omega_{\alpha}(x) = \{x(x-1)(x-t)\}^{-\alpha}dx
\]
according to \cite[Chap. IV]{Y97}, for $\alpha = \frac{k}{2k+1}$ and $\frac{2k-1}{4k}$ 
with $k \geq 2$. They satisfy differential equations 
$E(\frac{k}{2k+1}, \frac{k-1}{2k+1}, \frac{2k}{2k+1})$ and
$E(\frac{2k-1}{4k}, \frac{2k-3}{4k}, \frac{2k-1}{2k})$
with monodromy groups $\Delta(n,n,n)$, $n=2k+1$ and $2k$ respectively.
Let us consider decompositions
\begin{align*}
\PP^1(\CC) = \HH_+ \cup \PP^1(\RR) \cup \HH_-, \qquad
\PP^1(\RR) = I_0 \cup I_1 \cup I_2 \cup I_3, \\
I_0 = (-\infty, 0), \quad I_1 = (0,t), \quad I_2 = (t,1), \quad I_3 = (1, \infty)
\end{align*}
where $\HH_+$ and $\HH_-$ are the upper and lower half planes respectively, 
and $I_k$ are (oriented) real intervals. 
(As the initial position of $t$, we assume that $0 < t < 1$.) 
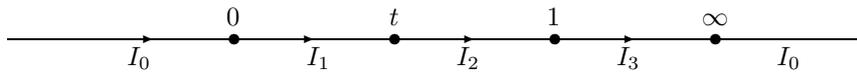
\begin{center}
\begin{picture}(300,40)(0,0)
\put(-10,20){\line(1,0){320}}
\multiput(35, 20)(60,0){4}{\vector(1,0){10}}
\multiput(75, 20)(60,0){4}{\circle*{4}}
\put(72, 25){\text{$0$}}
\put(133, 25){\text{$t$}}
\put(192, 25){\text{$1$}}
\put(250, 25){\text{$\infty$}}
\put(35, 10){\text{$I_0$}}
\put(102, 10){\text{$I_1$}}
\put(158, 10){\text{$I_2$}}
\put(218, 10){\text{$I_3$}}
\put(278, 10){\text{$I_0$}}
 \end{picture}
\end{center}
Modifying boundaries $\partial \HH_{+}$ and $\partial \HH_{-}$ to avoid 
$0, t, 1$ and $\infty$ as follows, we fix a branch of $\Omega_{\alpha}(x)$ on a 
simply connected domain $\HH_-$ and define integrals 
$\DS u_k(t) = \int_{I_k} \Omega_{\alpha}(x)$ by this branch.  
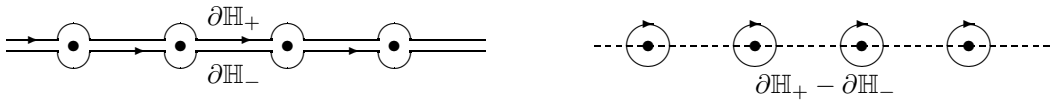
\begin{center}
\begin{picture}(440,40)(0,0)
\multiput(25,20)(40,0){4}{\circle*{4}} 
\multiput(240,20)(40,0){4}{\circle*{4}}
\multiput(25,22)(40,0){4}{\oval(12,13)[t]}
\multiput(25,18)(40,0){4}{\oval(12,13)[b]}
\put(0,22){\line(1,0){19}} \put(0,18){\line(1,0){19}}
\multiput(31,22)(40,0){4}{\line(1,0){28}} \multiput(31,18)(40,0){4}{\line(1,0){28}} 
\put(10, 22){\vector(1,0){2}} \put(50, 18){\vector(1,0){2}}
\put(90, 22){\vector(1,0){2}} \put(130, 18){\vector(1,0){2}}
\multiput(240, 28)(40,0){4}{\vector(1,0){2}}
\put(75, 28){\text{$\partial \HH_+$}} \put(75, 6){\text{$\partial \HH_-$}}
\put(280, 2){\text{$\partial \HH_+ - \partial \HH_-$}}
\multiput(220,20)(4,0){43}{\line(1,0){2}}
\multiput(240,20)(40,0){4}{\circle{16}}
\end{picture}
\end{center}
By the Cauchy integral theorem, they satisfy 
\begin{align*}
0 &= \int_{\partial \HH_-} \Omega_{\alpha}(x) 
 = u_0(t) + u_1(t) + u_2(t) + u_3(t),  
\\
0 &= \int_{\partial \HH_+} \Omega_{\alpha}(x)
   = u_0(t) + c u_1(t) + c^2 u_2(t) + c^3 u_3(t), \quad 
c = \exp(2\pi i \alpha)
\end{align*}
since $\Omega_{\alpha}(x)$ is multiplied by $\exp(2\pi i \alpha)$ 
if $x$ travels around $0, t$ or $1$ in clockwise direction. Hence we have
\begin{align*}
  u_2(t) = -\frac{1}{1+c} \{u_1(t) + (1+c+c^2) u_3(t) \}.
\end{align*}
\subsection{Monodromy}
Now let $\delta_0$ and $\delta_1$ be paths to make a half turn around $0$ and 
$1$ respectively in counter clockwise direction, starting from the initial 
point of $t$. 
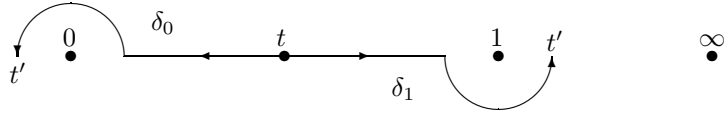
\begin{center}
\begin{picture}(420,60)(0,0)
\multiput(100,30)(80,0){4}{\circle*{4}}
\put(97,34){\text{$0$}}
\put(177,34){\text{$t$}}
\put(257,34){\text{$1$}}
\put(335,34){\text{$\infty$}}
\put(77,20){\text{$t'$}}
\put(278,32){\text{$t'$}}
\put(180,30){\line(1,0){60}}
\put(180,30){\line(-1,0){60}}
\put(100,30){\oval(40,40)[t]}
\put(260,30){\oval(40,40)[b]}
\put(80,32){\vector(0,-1){2}} \put(280, 28){\vector(0,1){2}}
\put(210,30){\vector(1,0){2}} \put(150, 30){\vector(-1,0){2}}
\put(130,40){\text{$\delta_0$}} \put(220,15){\text{$\delta_1$}}
\end{picture}
\end{center}
Corresponding analytic continuations are represented by connection matrices 
$h_1$ and $h_1$:
\begin{align*}
\delta_0 &:
   \begin{bmatrix} u_1(t) \\ u_3(t) \end{bmatrix} 
\dashrightarrow
   \begin{bmatrix} -c^{-1} u_1(t') \\ u_3(t') \end{bmatrix}
= h_0 \begin{bmatrix} u_1(t') \\ u_3(t') \end{bmatrix}, 
\qquad
h_0 = \begin{bmatrix} -c^{-1} & 0 \\ 0 & 1 \end{bmatrix} \\
\delta_1 &: 
  \begin{bmatrix} u_1(t) \\ u_3(t) \end{bmatrix} \dashrightarrow
   \begin{bmatrix} u_1(t') + u_2(t') \\ c^{-1} u_2(t') + u_3(t') \end{bmatrix}
= h_1 \begin{bmatrix} u_1(t') \\ u_3(t') \end{bmatrix}, \qquad
h_1 = \begin{bmatrix} \frac{c}{c+1} &-\frac{c^2+c+1}{c+1} \\ 
 -\frac{1}{c^2+c} & -\frac{1}{c^2+c} \end{bmatrix}
\end{align*}
where $u_1(t'), \dots, u_4(t')$ are integrals over oriented intervals 
$I_1', \cdots, I_4'$ defined for new configrations $-\infty < t'<0<1<\infty$ 
and $-\infty< 0 < 1<t'<\infty$. 
The monodromy group $\mathbf{Mon}$ is generated by 
\[
g_0 = h_0^2 = \begin{bmatrix} c^{-2} & 0 \\ 0 & 1 \end{bmatrix},
\qquad
g_1 = h_1^2 =
\begin{bmatrix}
 \frac{c^2+1}{c^2+c} & \frac{1-c^3}{c^2+c} \\
 \frac{1-c}{c^3+c^2} & \frac{c^2+1}{c^3+c^2} 
\end{bmatrix}. 
\]
\subsection{Hermitian form and period domain}
It is known that there exists a unique monodromy-invarinat Hermitian form up 
to constant (see e.g. \cite{B07} and \cite{Y97}). 
In fact, we can easily check that $h_0$ and $h_1$ belong to a unitary group
\[
 U_H = \{ g \in \mathrm{GL}_2(\CC) \ | \ {}^t \bar{g} H g = H \}, \qquad
 H = \begin{bmatrix} 1 & 0 \\ 0 & 1 + c + c^{-1} \end{bmatrix},
\]
and hence $\mathbf{Mon} \subset U_H$. The value of $1 + c + c^{-1}$ is negative 
for $c = \exp(2 \pi i \alpha)$ with  
$\alpha = \frac{k}{2k+1}$ and $\frac{2k-1}{4k}$ \ $(k \geq 2)$, and $H$ is 
indefinite. Therefore two domains 
\[
\DD_H^{\pm} = \{ u \in \CC^2 \ | \ \pm {}^t \bar{u} H u  < 0 \} / \CC^{\times} 
\subset \PP^1(\CC).
\]
are disks, and $U_H$ acts on each domain.  
Now the image of the Schwarz map
\[
 \mathfrak{s} : \CC-\{0, 1\} \longrightarrow \PP^1(\CC), \quad t \mapsto [u_1(t): u_3(t)]
\]
is contained in either $\DD_H^+$ or $\DD_H^-$, which is tessellated by Schwarz triangles. 
Since we have
\[
\mathfrak{s}(0) = \lim_{t \rightarrow 0}[u_1(u):u_3(t)] = [0:u_3(0)] \in \DD_H^+,
\]
we see that $\DD_H^+ / \mathbf{Mon} \cong \PP^1(\CC)$ and 
$\DD_H^+ / [\mathbf{Mon}, \mathbf{Mon}] \cong \mathcal{F}_n$, where 
$\mathcal{F}_n$ is the Fermat curve of degree $n$ with 
$n = 2k+1$ (resp. $2k$) if $\alpha = \frac{k}{2k+1}$ (resp. $\frac{2k-1}{4k}$).
\subsection{Remark}
(1) Putting $\zeta_d = \exp(2 \pi i/ d)$, we have
\[
1+c+c^{-1} =
\begin{cases}
  1 + (\z_{2k+1})^k + (\z_{2k+1})^{k+1} \qquad (n=2k+1) \\
  1 + (\z_{4k})^{2k-1} + (\z_{4k})^{2k+1} \qquad (n=2k)
\end{cases}
\]
(2) In the case of $n = 2k+1$, we have
\[
g_0 = \begin{bmatrix} \z & 0 \\ 0 & 1 \end{bmatrix}, \quad
g_1 = \frac{1}{1 + \z^k}
\begin{bmatrix}
  \z^k + \z^{k+1} & \z^{k+1} - \z^{2k} \\ \z - \z^{k+1} & 1 + \z
\end{bmatrix} 
\]
where $\z = \z_{2k+1}$. Since $1/(1 + \z^k) = -(\z + \z^2 + \dots + \z^k)$ and 
$\det g_1 = \z$, the monodromy group 
$\mathbf{Mon}$ is a subgroup of $U_H \cap \mathrm{GL}_2(\ZZ[\z])$.  
\\
(3) In the case of $n = 2k$, we have
\[
g_0 = \begin{bmatrix} \z^2 & 0 \\ 0 & 1 \end{bmatrix}, \quad
g_1 = \frac{1}{1 + \z^{2k-1}}
\begin{bmatrix}
  \z^{2k+1} + \z^{2k-1} & \z^{2k+1} - \z^{4k-2} \\ \z^2 - \z^{2k+1} & 1 + \z^2
\end{bmatrix} 
\]
where $\z = \z_{4k}$. Note that the cyclotomic polynomial $\Phi_{4k}(x)$ satisfies 
$\Phi_{4k}(1) = 1$ if $4k \ne 2^m$. In this case, $1-\z$ is a unit in 
$\ZZ[\z]$, and so is $1/(1 + \z^{2k-1}) = \z / (\z - 1)$. Hence  
$\mathbf{Mon}$ is a subgroup of $U_H \cap \mathrm{GL}_2(\ZZ[\z])$ if $4k \ne 2^m$. 
\subsection{Fermat curve as a Shimura variety}
A triangle group $\Delta(n,n,n)$ is arithmetic for 
\[
n \in FT = \{4, 5, 6, 7, 8, 9, 12, 15 \}, 
\]
and the Fermat curve $\mathcal{F}_n$ is a Shimura curve. Let us see corresponding 
families of hypergeometric curves
\[
X_t \ : \ y^m = x(x-1)(x-t)
\]
for these case. By the Riemann-Hurwitz formula, the genus of $X_t$ is $g = m-1$ 
if $3 \nmid m$, and $g = m-2$ if $3 \mid m$.
Let $\rho$ be the covering automorphism $(x,y) \rightarrow (x, \z_m y)$ where 
$\z_m = \exp(2 \pi i / m)$. By this action, we can decompose $H^1(X_t,\QQ)$ into
irreducible representations of $\rho$, and $H^1(X_t,\CC)$ into eigenspaces of $\rho$. 
Let $V(\lambda)$ be the $\lambda$-eigenspace of $\rho$. 
If $m$ is not prime, the covering $X_t \rightarrow \PP^1$ has intermediate curves 
$Y_t$, and the pullback of $H^1(Y_t, \CC)$ consists of $V(\z_m^k)$ such that $(m,k) \ne 1$. 
Conversely, such $V(\z_m^k)$ descends to a quotient curve. 
From explicit basis of $H^{1,0}(X_t)$, we see that the Prym part
\[
 H^1_{Prym}(X_t,\QQ) = [ \oplus_{(k,m)=1} V(\z_m^k) ] \cap H^1(X_t, \QQ) 
\]  
has a Hodge structure of type
\[
  H^1_{Prym}(X_t,\CC) 
= \underbrace{V(\lambda_1) \oplus \dots \oplus 
V(\lambda_{d-1})}_{\text{contained in} \ H^{1,0}}
\oplus \underbrace{V(\lambda_d)}_{\text{split}}
\oplus \underbrace{V(\lambda_{d+1})}_{\text{split}}
\underbrace{\oplus V(\lambda_{d+2}) \oplus \dots \oplus 
V(\lambda_{2d})}_{\text{contained in} \ H^{0,1}}
\] 
where $2d = [\QQ(\z_m):\QQ]$, $\lambda_1, \dots, \lambda_{2d}$ are primitive roots 
of unity $\z_m, \cdots, \z_m^{m-1}$ such that $\bar{\lambda_i} = \lambda_{2d+1-i}$ and 
$\dim V(\lambda_i) =2$ for $i = 1, \dots, 2d$. Therefore the Hodge strucure on 
$H^1_{Prym}(X_t,\QQ)$ with the action of $\rho$ is determined by a decomposition 
$V(\lambda_d) = V(\lambda_d)^{1,0} \oplus V(\lambda_d)^{0,1}$ (the docomposition of 
$V(\lambda_{d+1})$ is automatically determined as the complex conjugate 
of $V(\lambda_d)$, and vice versa), that is, determined by periods of
$\Omega_{\alpha}(x) \in V(\lambda_d)^{1,0}$.
\[
\begin{array}{|c||c|c|c|c|} 
\hline & & & & x^adx/y^b \ \text{with the following} \ $(a,b)$\\
  \Delta(n,n,n) & m & g & [\QQ(\z_m):\QQ] & 
\text{give a basis of} \ H^{1,0}(X_t)_{Prym} \\ \hline
  (4,4,4) & 8 & 7 & 4 & (0,3), (0,5), (0,7), (1,7) \\ \hline
  (5,5,5) & 5 & 4 & 4 & (0,2), (0,3), (0,4), (1,4) \\ \hline 
  (6,6,6) & 12 & 10 & 4 & (0,5), (0,7), (0,11), (1,11) \\ \hline
  (7,7,7) & 7 & 6 & 6 & (0,3), (0,4), (0,5), (1,5), (0,6),(1,6) \\ \hline
  (8,8,8) & 16 & 15 & 8 & (0,7), (0,9), (0,11), (1,11), \\
          & & & & (0,13),(1,13), (0,15),(1,15) \\ \hline
  (9,9,9) & 9 & 7 & 6 & (0,4), (0,5), (0,7), (1,7),(0,8),(1,8) \\ \hline
  (12,12,12) & 24 & 22 & 8 & (0,11), (0,13), (0,17), (1,17), \\ 
          & & & & (0,19),(1,19), (0,23),(1,23) \\ \hline
  (15,15,15) & 15 & 13 & 8 & (0,7), (0,8), (0,11), (1,11), \\
          & & & & (0,13),(1,13), (0,14),(1,14) \\ \hline
\end{array}
\]
\subsection{}
In the cases $n = 5$ and $7$, the monodromy group has a nice representatoin.
Put
\[
\Gamma = U_H \cap \mathrm{GL}_2(\ZZ[\z_n]), \qquad
\Gamma(\mathfrak{m}) = \{ g \in \Gamma \ | \ g \equiv 1 \mod \mathfrak{m}\}.
\]
The arithmetic quotient $\DD_H^+ / \Gamma$ is the moduli space of 
Jacobians of curves $y^n = x^3 + ax +b \ (n=5,7)$ as a PEL-family (see \cite{Sm64}).
Therefore we have the following diagram
\begin{align*}
  \begin{CD}
\DD_H^+ / \mathbf{Mon} @>>> \PP^1 = 
\overline{\{ \text{ordered distinct $3+1$ points} \ (0,1,t,\infty)\}}\\
     @VVV @VVV \\
\DD_H^+ / \Gamma @>>> \PP^1 / S_3 = \overline{\{\text{unordered distinct 3 points in} 
\ \CC \} / \sim}
  \end{CD}
\end{align*}
where horizontal arrow are isomorphisms, and $\sim$ is the equivalence relation 
by affine transformations. From this fact, we see that $\Gamma / \mathbf{Mon}$ is 
isomorphic to $S_3$ up to the center.
\subsection{Remark} 
For $n=5$, the Hermitian form $H$ is same with one given in \cite{Sm64}:
\[
 H = \begin{bmatrix} 1 & 0 \\ 0 & 1 + \z_5^2 + \z_5^3 \end{bmatrix}
= \begin{bmatrix} 1 & 0 \\ 0 & (1 - \sqrt{5})/2 \end{bmatrix}.
\]
For $n=7$, the Hermitian form given in \cite{Sm64} is
\[
 S = \begin{bmatrix} 1 & 0 \\ 0 & -\frac{\sin(3 \pi/7)}{\sin(2 \pi/7)} \end{bmatrix}
= \begin{bmatrix} 1 & 0 \\ 0 & -(\z_7 + \z_7^6) \end{bmatrix}
= {}^t \bar{A}HA, \qquad A = \begin{bmatrix} 1 & 0 \\ 0 & \z_7 + \z_7^6 \end{bmatrix}
\in \mathrm{GL}_2(\ZZ[\z_7]).
\]
\subsection{Proposition (\cite{YY84} for $n=5$)} \label{Prop-group}
Let us denote the image of 
$G \subset \mathrm{GL}_2(\ZZ[\z_n])$ in $P\mathrm{GL}_2(\ZZ[\z_n])$ by $\overline{G}$. 
For $n=5$ and $7$, 
\\
(1) the projective modular group $\overline{\Gamma}$ is projectively generated 
by $h_0$ and $h_1$,
\\
(2) we have 
\[
 \overline{\mathbf{Mon}} = \overline{\Gamma(1-\z_n)}, \qquad  
\overline{[\mathbf{Mon},\mathbf{Mon}]} = \overline{\Gamma((1-\z_n)^2)}
\] 
 as automorphisms of $\DD_H^+$.
\vskip0.3cm \noindent
{\it Proof}.
We show these facts only for $n=7$, but the case $n=5$ is shown by the same way 
(also see \cite{YY84} and \cite{K03} for $n=5$). The quotient group 
$\Gamma / \Gamma(1-\z_7)$ is isomorphic to a subgroup of the finite orthogonal group
\[
 \mathrm{O}(Q,\FF_7) = \{ g \in \mathrm{GL}_2(\FF_7) \ | \ {}^tgQg = Q \}, 
\qquad \FF_7 = \ZZ[\z_7] / (1-\z_7), \quad
Q = \begin{bmatrix} 1 & 0 \\ 0 & 3 \end{bmatrix}.
\]
The group $\mathrm{O}(Q,\FF_7)$ is isomorphic to $S_3 \times \{  \pm 1 \}$, 
since elements of $\mathrm{O}(Q,\FF_7) / \{ \pm 1 \}$ are
\begin{align*}
\text{order} \ 2: \quad 
 \begin{bmatrix}-1 & 0 \\ 0 & 1\end{bmatrix}, \quad
 \begin{bmatrix}3 & 2 \\ 3 & 4\end{bmatrix}, \quad
 \begin{bmatrix}4 & 2 \\ 3 & 3\end{bmatrix},
\qquad
\text{order} \ 3: \quad
\begin{bmatrix}3 & 2 \\ 4 & 3\end{bmatrix}, \quad
\begin{bmatrix}3 & 5 \\ 3 & 3\end{bmatrix}.
\end{align*}
Since we have
\begin{align*}
  h_0 \equiv \begin{bmatrix} -1 & 0 \\ 0 & 1 \end{bmatrix}, \quad
  h_1 \equiv \begin{bmatrix} 4 & 2 \\ 3 & 3 \end{bmatrix} \quad \mod 1-\zeta_7,
\end{align*}
the group $\Gamma / \Gamma(1-\z_7)$ is generated by $h_0, \ h_1$ and $\pm1$, 
and isomorphic to $S_3 \times \{  \pm 1 \}$. 
Therefore $\overline{\mathbf{Mon}}$ coincides with $\overline{\Gamma(1-\z_7)}$ 
since we have $\mathbf{Mon} \subset \Gamma(1-\z_7)$ and 
$\overline{\Gamma} / \overline{\mathbf{Mon}} = S_3$.  
Note that $\mathbf{Mon}$ is generated by $h_0^2$ and $h_1^2$, and hence 
$\overline{\Gamma}$ is generated by $h_0$ and $h_1$. A homomorphism
\begin{align*}
  \nu : \Gamma(1-\z_7) \longrightarrow \mathrm{M}_2(\FF_7), 
\qquad
 \nu(g) = \frac{1}{1-\z_7} (g - 1) \mod 1-\z_7
\end{align*}
has the kernel $\Gamma((1-\z_7)^2)$, and the image is 
generated by
\begin{align*}
  \nu(g_0) = \begin{bmatrix} -1 & 0 \\ 0 & 0 \end{bmatrix}, \quad
    \nu(g_1) = \begin{bmatrix} 5 & 1 \\ 5 & 1 \end{bmatrix}.
\end{align*}
Therefore we have $\Gamma(1-\z_7) / \Gamma((1-\z_7)^2) \cong (\ZZ / 7 \ZZ)^2$. 
Since we have
\[
[\Gamma(1-\z_7), \Gamma(1-\z_7)] \subset \Gamma((1-\z_7)^2), \qquad
\mathbf{Mon} / [\mathbf{Mon},\mathbf{Mon}] \cong (\ZZ / 7 \ZZ)^2,
\] 
we conclude that $\overline{[\mathbf{Mon},\mathbf{Mon}]} = 
\overline{\Gamma((1-\z_n)^2)}$.
\qed
\section{Heptagonal Curves}
\subsection{}
From now, we concentrate in the case $n = 7$, that is, a 1-dimensional family 
of algebraic curves 
\[
 X_t \ : \ y^7 = x(x-1)(x-t).
\] 
We denote $\z_7 = \exp(2 \pi i /7)$ simply by $\z$. As a Riemann surface, $X_t$ is obtained 
by glueing seven sheets $\Sigma_1, \cdots, \Sigma_7$, each of which is a copy of $\PP^1$ 
with cuts (see Figure 1) and satisfying $\rho(\Sigma_i) = \Sigma_{i+1}$ where 
indices are consdered modulo $7$. Let $\mathfrak{i}_i(x_1, x_2)$ be an oriented real interval 
from $x_1$ to $x_2$ on $\Sigma_i$. We define $1$-cycles
\begin{align*}
\gamma_1 = \mathfrak{i}_1(0,t) + \mathfrak{i}_2(t,0) = (1-\rho)\mathfrak{i}_1(0,t), \qquad 
\gamma_2 = \mathfrak{i}_1(t,1) + \mathfrak{i}_2(1,t) = (1-\rho)\mathfrak{i}_1(t,1), \\
\gamma_3= \mathfrak{i}_1(1,\infty) + \mathfrak{i}_2(\infty,1) = (1-\rho)\mathfrak{i}_1(1,\infty). 
\end{align*}
For computation of intersection 
numbers, we use deformations of $\gamma_1$ and $\gamma_3$ as in Figure 1. 
\begin{figure}[htbp] \begin{center}
\setlength{\unitlength}{1mm}
\begin{picture}(135,40)
\thinlines
\multiput(10,5)(0,1){7}{\line(0,1){0.5}} \multiput(55,5)(0,1){7}{\line(0,1){0.5}}
\multiput(90,5)(0,1){7}{\line(0,1){0.5}} \multiput(135,5)(0,1){7}{\line(0,1){0.5}} 
\multiput(10,12)(1,0){125}{\line(1,0){0.5}}
\multiput(10,20)(0,1){7}{\line(0,1){0.5}} \multiput(55,20)(0,1){7}{\line(0,1){0.5}}
\multiput(90,20)(0,1){7}{\line(0,1){0.5}} \multiput(135,20)(0,1){7}{\line(0,1){0.5}}
\multiput(10,27)(1,0){125}{\line(1,0){0.5}}
\multiput(10,35)(0,1){7}{\line(0,1){0.5}} \multiput(55,35)(0,1){7}{\line(0,1){0.5}}
\multiput(90,35)(0,1){7}{\line(0,1){0.5}} \multiput(135,35)(0,1){7}{\line(0,1){0.5}}
\multiput(10,42)(1,0){125}{\line(1,0){0.5}}
\put(0,5){$\Sigma_1$} \put(0,20){$\Sigma_2$} \put(0,35){$\Sigma_5$}
\put(46,7){to $\Sigma_2$} \put(11,22){to $\Sigma_1$}
\put(126,7){to $\Sigma_5$} \put(91,22){to $\Sigma_1$}
\multiput(10,5)(0,15){3}{\circle*{1.2}} \put(9,0){$0$} 
\multiput(55,5)(0,15){3}{\circle*{1.2}} \put(54,0){$t$} 
\multiput(90,5)(0,15){3}{\circle*{1.2}} \put(89,0){$1$} 
\multiput(135,5)(0,15){3}{\circle*{1.2}} \put(133,0){$\infty$}
\thicklines
\put(10,5){\oval(5,5)[tr]}
\put(12.5,5){\line(1,0){40}} \put(15,5){\vector(1,0){5}} 
\put(55,5){\oval(5,5)[b]} \put(55,5){\oval(5,5)[tr]}
\put(10,20){\oval(5,5)[l]} \put(10,20){\oval(5,5)[br]}
\put(12.5,20){\line(1,0){40}} \put(40,20){\vector(-1,0){5}}
\put(55,20){\oval(5,5)[tl]} 
 \put(30,35){$\gamma_1$} 
\thicklines
\put(90,5){\oval(5,5)[tr]}
\put(92.5,5){\line(1,0){40}} \put(95,5){\vector(1,0){5}} 
\put(135,5){\oval(5,5)[b]} \put(135,5){\oval(5,5)[tr]}
\put(90,20){\oval(5,5)[l]} \put(90,20){\oval(5,5)[br]} 
\put(92.5,20){\line(1,0){40}} \put(125,20){\vector(-1,0){5}} 
\put(135,20){\oval(5,5)[tl]}
\put(135,35){\oval(5,5)[l]} \put(135,35){\oval(5,5)[r]} 
\put(110,35){$\gamma_3$}
\put(132.5,35){\vector(0,-1){1}}
\end{picture} \end{center} \caption{} \label{Fgcycle}
\end{figure}
\vskip1mm
Let $\mathbf{Int}_k$ be the intersection matrix 
$[\rho^i(\gamma_k) \cdot \rho^j(\gamma_k)]_{0 \leq i,j \leq 5}$. We have 
\begin{align*}
\mathbf{Int}_1 =
\begin{bmatrix}
 0 & 1 & 0 & 0 & 0 & 0 \\
 -1 & 0 & 1 & 0 & 0 & 0 \\
 0 & -1 & 0 & 1 & 0 & 0 \\
 0 & 0 & -1 & 0 & 1 & 0 \\
 0 & 0 & 0 & -1 & 0 & 1 \\
 0 & 0 & 0 & 0 & -1 & 0 
\end{bmatrix}
\qquad 
\mathbf{Int}_3 =
\begin{bmatrix}
 0 & 1 & 0 & -1 & 1 & 0 \\
 -1 & 0 & 1 & 0 & -1 & 1 \\
 0 & -1 & 0 & 1 & 0 & -1 \\
 1 & 0 & -1 & 0 & 1 & 0 \\
 -1 & 1 & 0 & -1 & 0 & 1 \\
 0 & -1 & 1 & 0 & -1 & 0 
\end{bmatrix}
\end{align*}
and $\det \mathbf{Int}_1 = \det \mathbf{Int}_3 = 1$. Since 
$\rho^i(\gamma_1) \cdot \rho^j(\gamma_3) = 0$, the intersecion matrix of twelve $1$-cycles 
$\gamma_1, \rho(\gamma_1), \cdots, \rho^5(\gamma_1)$ and 
$\gamma_3, \rho(\gamma_3), \cdots, \rho^5(\gamma_3)$ is unimodular, and they 
form a basis of $H_1(X_t, \ZZ)$. 
Hence $\{ \gamma_1, \ \gamma_3 \}$ gives a basis of $H_1(X_t, \ZZ) \cong \ZZ[\rho]^2$ as 
a $\ZZ[\rho]$-module.
\\ \indent
Similarly we have $H^1(X_t, \ZZ) \cong \ZZ[\rho]^2$ and
the decomposition of $H^1(X_t, \CC) \cong \ZZ[\rho]^2 \otimes \CC$ 
into eigenspaces of $\rho$:
\[
 H^1(X_t, \CC) = V(\z) \oplus V(\z^2) \oplus \cdots \oplus V(\z^6), \qquad
\dim V(\z^k)=2.
\]
Let $P_0, \ P_1, P_t$ and $P_{\infty}$ be four ramification points of $X_t$ over $0, 1, t$ 
and $\infty$. We denote the divisor of a rational function (or a rational $1$-form) 
$f$ by $\mathrm{div}(f)$. 
Then we see that
\begin{align*}
\mathrm{div}(x) = 7P_0 - 7P_{\infty}, \qquad
\mathrm{div}(y) = P_0 + P_1 + P_t - 3P_{\infty}, \\
\mathrm{div}(dx) = 6(P_0 + P_1 + P_t) - 8P_{\infty},
\end{align*}
and holomorphic $1$-forms
\begin{align*}
\omega_1 = \frac{dx}{y^3}, \quad \omega_2 = \frac{dx}{y^4}, \quad 
\omega_3 = \frac{dx}{y^5}, \quad \omega_4 = \frac{xdx}{y^5}, \quad 
\omega_5 = \frac{dx}{y^6}, \quad \omega_6 = \frac{xdx}{y^6}
\end{align*}
on $X_t$ give a basis of $H^{1,0}(X_t)$. 
\subsection{Remark}
As stated in the previous section, we have
\[
 V(\z) \oplus V(\z^2) \subset H^{1,0}(X_t), \qquad 
 V(\z^5) \oplus V(\z^6) \subset H^{0,1}(X_t)
\] 
and the Hodge structure on $H^1(X_t, \ZZ)$ is determined by a decomposition of 
$V(\z^4)$. 
\subsection{Period Matrix}
The following $1$-cycles
\begin{align*}
\begin{array}{lll}
 B_1 = \gamma_1, & 
 B_2 = (1 + \rho^2) (\gamma_1), &
 B_3 = (1 + \rho^2 + \rho^4) (\gamma_1), 
\\
 A_1 = \rho(\gamma_1), & 
 A_2 = \rho^3(\gamma_1), &
 A_3 =\rho^5(\gamma_1),
\\
 B_4 = \rho^5(\gamma_3), &
 B_5 =  \rho^3(\gamma_3), &
 B_6 = (1+ \rho  - \rho^4 - \rho^5) (\gamma_3), 
\\
 A_4 = (1 + \rho^2) (\gamma_3), &
 A_5 = (-\rho + \rho^4 + \rho^5) (\gamma_3), &
 A_6 = (1 + \rho + \rho^2) (\gamma_3).
\end{array}
\end{align*}
give a symplectic basis of $H_1(X_t, \ZZ)$ such that 
\[
A_i \cdot A_j = 0, \quad B_i \cdot B_j = 0,
\quad B_i \cdot A_j = \delta_{ij}. 
\]
The associated period matrix is
\begin{align*}
 \Pi_A = \begin{bmatrix} \int_{A_i} \omega_j \end{bmatrix}
 = \begin{bmatrix}
 \int_{\gamma_1} \vec{\omega} R \\ \int_{\gamma_1} \vec{\omega} R^3 \\ 
\int_{\gamma_1} \vec{\omega} R^5 \\ \int_{\gamma_3} \vec{\omega} (I + R^2) \\ 
\int_{\gamma_3} \vec{\omega}(-R + R^4 + R^5) \\ \int_{\gamma_3} \vec{\omega} (I + R + R^2)
 \end{bmatrix},
\qquad
 \Pi_B = \begin{bmatrix} \int_{B_i} \omega_j \end{bmatrix}
 = \begin{bmatrix}
 \int_{\gamma_1} \vec{\omega} \\ \int_{\gamma_1} \vec{\omega} (I + R^2) \\ 
\int_{\gamma_1} \vec{\omega} (I + R^2 + R^4) \\ \int_{\gamma_3} \vec{\omega} R^5 \\ 
\int_{\gamma_3} \vec{\omega} R^3 \\ \int_{\gamma_3} \vec{\omega} (I + R - R^4 - R^5)
 \end{bmatrix}
\end{align*}
where $\vec{\omega} = (\omega_1, \dots, \omega_6)$ and 
$R = \mathrm{diag}(\z^4, \z^3, \z^2, \z^2, \z, \z)$. 
The normalized period matrix $\tau = \Pi_A \Pi_B^{-1}$ belongs to the Siegel 
upper half space $\HH_6$, consisting of symmetric matrices of degree $6$ whose imaginary 
part is positive definite. The symplectic group 
\[
Sp_{12}(\ZZ) = \{ \gamma \in \mathrm{GL}_{12}(\ZZ) \ | \ {}^t \gamma J \gamma = J \}, \quad 
J = \begin{bmatrix} 0 & \I_6 \\ -\I_6 & 0\end{bmatrix},
\]
acts on $\HH_6$ by $\begin{bmatrix} a & b \\ c & d \end{bmatrix} \cdot \tau 
= (a \tau + b)(c \tau + d)^{-1}$, and $\mathcal{A}_6 = \HH_6 / Sp_{12}(\ZZ)$ is 
the moduli space of principally polarized abelian varieties (p.p.a.v.) of dimension $6$.
\subsection{Remark}
For a suitable choice of a branch of $\Omega_{\alpha}(x)$ in the previous section, 
we have 
\[
 \int_{\gamma_k} \omega_1 = (1-\z^4) u_k(t) \qquad (k=1,2,3).
\]
Since we use $u_k$ for projective coordinates mainly, hereafter we denote 
$\int_{\gamma_k} \omega_1$ by $u_k$ for simplicity.
\subsection{Symplectic representation}
Let $M \in Sp_{12}(\ZZ)$ be the symplectic representation of $\rho$ with respect to 
the above basis:
\[
 (\rho(A_1), \dots, \rho(A_6), \rho(B_1), \dots, \rho(B_6))
= (A_1, \dots, A_6, B_1, \dots, B_6)^tM.
\]
Explicit form of $M$ is given in Appendix. By definition, we have 
$M \begin{bmatrix} \Pi_A \\ \Pi_B \end{bmatrix} =
 \begin{bmatrix} \Pi_A \\ \Pi_B \end{bmatrix} R$. 
Therefore $\Pi_A \Pi_B^{-1}$ belongs to a domain 
$\HH_6^M = \{ \tau \in \HH_6 \ | \ M \cdot \tau = \tau \}$, which parametrizes 
p.p.a.v of dimension $6$ with an automorphism $M$
(see section 5 in \cite{vG92}). We know that this domain is $1$-dimensional, 
and hence isomorphic to $\DD_H^+$ (\cite{BL92}, Chap. 9 and \cite{Sm64}). 
The centralizer of $M$ in $Sp_{12}(\ZZ)$
\[
 Sp_{12}^M(\ZZ) = \{ g \in Sp_{12}(\ZZ) \ | \ g M = M g \}.
\]
acts on the domain $\HH_6^M$.
\subsection{Proposition} \label{modular-embedding}
There exist a group isomorphisms $\phi : \Gamma \rightarrow Sp_{12}^M(\ZZ)$ and 
an analytic isomorphism $\Phi : \DD_H^+ \rightarrow \HH_6^M$ such that 
$\Phi(gu) = \phi(g) \Phi(u)$, that give the following commutative diagram.
\begin{align*}
  \begin{CD}
\DD_H^+  @>\Phi>> \HH_6^M  \\
     @VVV @VVV \\
\DD_H^+ / \Gamma @>>> \HH_6^M / Sp_{12}^M(\ZZ)
  \end{CD}
\end{align*}
\vskip0.3cm \noindent
{\it Proof}.
 Now we have
\begin{align*}
\Pi_{A,1} = {}^t[\int_{A_1} \omega_1, \cdots, \int_{A_6} \omega_1] 
 ={}^t[\z^4 u_1,\ \z^5 u_1,\ \z^6 u_1,\ (1+\z)u_3,\ (\z^2 - \z^4 + \z^6)u_3,\ (1+\z+\z^4)u_3],
\\
\Pi_{B,1} = {}^t[\int_{B_1} \omega_1, \cdots, \int_{B_6} \omega_1]
 = {}^t[u_1, (1 + \z)u_1,\ (1+\z+\z^2)u_1,\ \z^6 u_3,\ \z^5 u_3,\ (1+\z^4 - \z^2 - \z^6)u_3].
\end{align*}
This correspondence 
$\begin{bmatrix} u_1 \\ u_3 \end{bmatrix} \mapsto 
\begin{bmatrix} \Pi_{A,1} \\ \Pi_{B,1} \end{bmatrix}$ 
define a linear map $\Phi_1 : \CC^2 \rightarrow \CC^{12}$. 
Since coefficients of $u_1$ (or $u_3$) in $\Pi_{A,1}$ and $\Pi_{B,1}$ give a 
$\ZZ$-basis of $\ZZ[\z]$, there exists a homomorphism 
$\phi : \mathrm{GL}_2(\ZZ[\z]) \rightarrow \mathrm{GL}_{12}(\ZZ)$ such that 
$\Phi_1(gu) = \phi(g) \Phi_1(u)$. Especialy, we have $\phi(\z^4 I_2) = M$ and the image of 
$\phi$ is the centralizer of $M$. We can easily check that the condition 
\[
|u_1|^2 + (1+\z^3 +\z^4)|u_3|^2 < 0
\] 
for $\DD_H^+$ is equivalent to Riemann's relation (\cite{M83})
\[
 \mathrm{Im} 
\left( \sum_{i=1}^6 \overline{\int_{B_i} \omega_1} \int_{A_i} \omega_1 \right) >0,
\]
and hence $\phi(\Gamma) = Sp_{12}^M(\ZZ)$. We give the map $\Phi$, which is compatible 
with $\Phi_1$, explicitly in Appendix.
\qed
\subsection{Remark} \label{Rem Riemann form}
Let us define a homomorphism
\begin{align*}
 \lambda : H_1(X_t, \ZZ) = \left< \gamma_1, \gamma_3 \right>_{\ZZ[\rho]} 
\longrightarrow \ZZ[\z]^2, \qquad
F_1(\rho)\gamma_1 + F_3(\rho) \gamma_3 \mapsto (F_1(\z^4), F_3(\z^4)). 
\end{align*}
By explicit computation, we see that the intersection form (which gives the polarization) 
on $H_1(X_t, \ZZ)$ is given by
\[
 E(x, y) = \frac{1}{7} Tr_{\QQ(\z)/\QQ}
((\z^3-\z^4) ^t \overline{\lambda(x)} H^{-1} \lambda(y)).
\]
\section{Schwarz inverse and theta function}
\subsection{Abel-Jacobi map}
For the normalized holomorphic $1$-forms
\[
 \vec{\xi} = (\xi_1, \dots, \xi_6 ) = (\omega_1, \dots, \omega_6) \Pi_B^{-1}
\]
with respect to $A_i$ and $B_i$ in the previous section, period integrals satify
\[
\tau = [\int_{A_i} \vec{\xi} \medspace]_{1 \leq i \leq 6} \in \HH_6^M, \qquad 
[\int_{B_i} \vec{\xi} \medspace]_{1 \leq i \leq 6} = I_6.
\]
Let $Div(X_t)$ be the group of divisors on $X_t$, and $J(X_t)$ be the Jacobian variety 
$\CC^6 / \ZZ^6 \tau + \ZZ^6$. The Abel-Jacobi map with the base point $P_{\infty}$ is 
\[
Div(X_t) \longrightarrow J(X_t), \qquad 
\sum m_i Q_i \mapsto \sum m_i \int_{P_{\infty}}^{Q_i} \vec{\xi} \mod \ZZ^6 \tau + \ZZ^6.
\]
We denote this homomorphism by $\overline{\fA}$, and a lift of $\overline{\fA}(D)$ 
by $\fA(D)$ (Hence $\fA : Div(X_t) \rightarrow \CC^6$ is a multi-valued map). 
As is well known, $\overline{\fA}$ factors through 
\[
Div(X_t) \longrightarrow Pic(X_t) = Div(X_t)/\{ \text{principal divisors} \}.
\]  
Since the base point is fixed by $\rho$, the map $\overline{\fA}$ is $\rho$-equivariant. 
Therefore the image of a $\rho$-invariant divisor belongs to the fixed points of $\rho$, 
that is, the $(1-\rho)$-tosion subgroup 
\[
 J(X_t)_{1-\rho} = \{ z \in J(X_t) \ | \ (1-\rho)z = 0\}.
\]
\subsection{Lemma} 
The $(1-\rho)$-tosion subgroup is 
\[
 J(X_t)_{1-\rho} = \{ \overline{\fA}(mP_0 + nP_1) \ | \ m,n \in \ZZ  \} \cong (\ZZ/ 7 \ZZ)^2
\] 
More explicitly, we have 
\[ 
  \fA(mP_0 + nP_1) \equiv a_{m,n}\tau + b_{m,n} \mod \ZZ^6 \tau + \ZZ^6
\] with
\begin{align*}
  a_{m,n} &= \frac{1}{7}(m, \, 2m,\, 3m,\, 2m + 3n,\, 2m + 3n,\, 0) \in \frac{1}{7}\ZZ^6,
\\
  b_{m,n} &= \frac{1}{7}(-m,\, -m,\, -m,\, 3m + n,\, 5m + 4n,\, m + 5n) \in \frac{1}{7}\ZZ^6.
\end{align*}
{\it Proof}.
It is obvious that $Ker (1-\rho) \cong (\ZZ[\z]/(1-\z))^2 \cong (\ZZ /7\ZZ)^2$.
Recall that
\begin{align*}
 \gamma_1 = (1-\rho)\mathfrak{i}_1(0,1), \quad
 \gamma_2 = (1-\rho)\mathfrak{i}_1(1,\infty), \quad
 \gamma_3 = (1-\rho)\mathfrak{i}_1(t,1).
\end{align*}
Computing intersection numbers, we see that
\[
\gamma_2 = A_1 + A_2 + A_3 + B_4 + B_5  =  \rho(\gamma_1) + \rho^3(\gamma_1) 
+ \rho^5(\gamma_1) + \rho^5(\gamma_3) + \rho^3(\gamma_3).
\]
Therefore we have
\begin{align*}
 \mathfrak{i}_1(0,t) 
&= \frac{1}{7}(6 + 5\rho + 4\rho^2 + 3\rho^3 + 2\rho^4 + \rho^5) \gamma_1
 = \frac{1}{7}(5A_1 + 3A_2 + A_3 + 2B_1 + 2B_2 + 2B_3), \\
 \mathfrak{i}_1(1,\infty) 
&= \frac{1}{7}(6 + 5\rho + 4\rho^2 + 3\rho^3 + 2\rho^4 + \rho^5) \gamma_3
 = \frac{1}{7}( -3A_4 + 4A_5 + 7A_6 - B_4 + 3B_5 + 2B_6), \\
 \mathfrak{i}_1(t,1) 
&= \frac{1}{7}(6 + 5\rho + 4\rho^2 + 3\rho^3 + 2\rho^4 + \rho^5) \gamma_2 
\\
 &= \frac{1}{7}(A_1 + 2A_2 + 3A_3 - B_1 - B_2 - B_3) 
 + \frac{1}{7}(A_4 + A_5 -7A_6 + 5B_4 -B_5 + 4B_6),
\end{align*}
namely,
\begin{align*}
\int^t_0 \vec{\xi} 
\equiv \frac{1}{7}(5,3,1,0,0,0)\tau + \frac{1}{7}(2,2,2,0,0,0), \quad
\int^{\infty}_1 \vec{\xi}
\equiv \frac{1}{7}(0,0,0,4,4,0)\tau + \frac{1}{7}(0,0,0,6,3,2), \\
\int^1_t \vec{\xi} 
\equiv \frac{1}{7}(1,2,3,1,1,0)\tau + \frac{1}{7}(6,6,6,5,6,4)
\mod \ZZ^6 + \tau \ZZ^6.
\end{align*}
As combinations of these integrals, we obtain explicit values of 
$\overline{\fA}(P_0)$ and $\overline{\fA}(P_1)$. 
\qed
\subsection{Theta function and Riemann constant}
Let us consider Riemann's theta function 
\[
  \vartheta(z,\tau) = \sum_{n \in \ZZ^6}\exp[\pi i n \tau {}^tn + 2 \pi i n{}^tz], 
\quad (z, \tau) \in \CC^6 \times \HH_6.
\]
The Abel-Jacobi map $\overline{\fA}$ induces a birational morphism from 
$\mathrm{Sym}^6 X_t$ to $J(X_t)$, and $W_{\fA}^5 = \overline{\fA}(\mathrm{Sym}^5 X_t)$ 
is a translation of the theta divisor
\[
  \Theta = \{ z \in J(X_t) \ | \ \vartheta(z) = 0 \}.
\]
More precisely, there exist a constan vector $\kappa \in \CC^6$ such that 
$\vartheta(e, \tau) = 0$ if and only if  
\begin{align*}
e \equiv \kappa - \fA(Q_1 + \cdots + Q_5) \mod \ZZ^6 \tau + \ZZ^6   
\end{align*}
for some $Q_1, \dots, Q_5 \in X_t$. 
The constant $\kappa$ (or its image $\overline{\kappa}$ in $J(X_t)$) is called 
the Riemann constant. It is the image of a half canonical class by $\fA$ 
(\cite{M83}, Chap. II, Appendix to \S 3), and depends only on a symplectic basis 
$A_i, B_i$ and the base point of $\fA$. Since 
$\mathrm{div}(\omega_5) = 10 P_{\infty}$, the image of the canonical class by 
$\overline{\fA}$ is $0$ and $\kappa$ must be a half period. Hence we have 
$\kappa = a \tau + b$ for some $a, b \in \frac{1}{2} \ZZ^6$. 
By the same argument as the proof of Lemma 5.4 in \cite{K03}, the corresponding theta 
characteristic $(a, b)$ is invariant under the action of $M$ on $\QQ^{12} / \ZZ^{12}$:
\[
M \cdot (a, b) = (a, b)M^{-1}
+ \frac{1}{2} ( \mathrm{diag}(C^tD), \mathrm{diag}(A^tB)), \qquad
M = \begin{bmatrix} A & B \\ C & D \end{bmatrix}.
\]
By explicit computation, we have 
\subsection{Lemma} \label{Lemma-torsion}
The $M$-invariant theta characteristics are $(a_{m,n} + a_0, b_{m,n} + b_0)$ with 
\[
a_0 = \frac{1}{2}(1,0,1,0,0,1), \quad b_0 = \frac{1}{2}(1,1,1,0,1,0).
\]
Especially, we have $\kappa \equiv a_0 \tau + b_0$. Since 
$\vartheta(-e) = \vartheta(e)$ and $\kappa$ is a half period, we have 
\[
 \overline{\kappa} - W_{\fA}^5 = \Theta = -\Theta  = \overline{\kappa} + W_{\fA}^5, 
\]
that is $W_{\fA}^5 = - W_{\fA}^5$.
\subsection{}
Let us consider $J(X)_{1-\rho} \cap W_{\fA}^5$. By definition, we hvae 
\[
\overline{\fA}(mP_0 + nP_1) \in W_{\fA}^5 = -W_{\fA}^5
\] 
for $0 \leq m, n \leq 6$ such that $m + n \leq 5$ or $(7-m) + (7-n) \leq 5$. 
The rest of $J(X)_{1-\rho}$ are $\overline{\fA}(mP_0 + nP_1)$ with the following $(m,n)$:
\begin{align*}
(1,5), \ (1,6), \ (2,4) \ (2,5), \ (2,6), \ (3,3), \ (3,4), \ (3,5), \ \\
(4,2), \ (4,3), \ (4,4), \ (5,1), \ (5,2), \ (5,3), \ (6,1), \ (6,2).
\end{align*}
Moreover we have the following reduction:
\begin{align*}
  (6P_0 + P_1) = (2P_1 + P_t + 4P_{\infty}) + \dv (\frac{x}{y}), \qquad
  (3P_0 + 3P_1) = (4P_t + 2P_{\infty}) + \dv (\frac{x(x-1)}{y^4}), \\
  (5P_0 + P_1) = (3P_1 + 2P_t + P_{\infty}) + \dv (\frac{x}{y^2}), \qquad
  (4P_0 + 3P_1) = (P_0 + 4P_t + 2P_{\infty}) + \dv (\frac{x(x-1)}{y^4}), 
\end{align*} 
that is, 
\[
\overline{\fA}(6P_0 + P_1), \quad \overline{\fA}(3P_0 + 3P_1), \quad
 \overline{\fA}(5P_0 + P_1), \quad \overline{\fA}(4P_0 + 3P_1) \in W_{\fA}^5.
\]
By the equality $W_{\fA}^5 = -W_{\fA}^5$ and symmetry for $P_0, P_1$,
 we see that $\overline{\fA}(mP_0 + nP_1) \in W_{\fA}^5$ if
\[
 (m,n) \ne (2,4),\ (2,5), \ (3,5), \ (4,2), \ (5,2), \ (5,3).
\]
The converse is also true:
\subsection{Lemma}
We hvae $\overline{\fA}(mP_0 + nP_1) \notin W_{\fA}^5$ for
\[
 (m,n) = (2,4),\ (2,5), \ (3,5), \ (4,2), \ (5,2), \ (5,3).
\]
{\it Proof}.
To prove this, note that
\[
 (5P_0 + 2P_1) = (4P_1 + 2P_t + P_{\infty}) + \dv (\frac{x}{y^2}), \qquad
  \therefore \ \overline{\fA}(5P_0 + 2P_1) = \overline{\fA}( 4P_1 + 2P_t)
\]
and 
\begin{align*}
  \overline{\fA}(3P_i + 5P_j) = - \overline{\fA}(4P_i + 2P_j), \qquad i,j \in \{0,1\}.
\end{align*}
By symmetry for $P_0, P_1$ and $P_t$, it suffices to prove that 
$\overline{\fA}(4P_0 + 2P_1) \notin W_{\fA}^5$. 
\\ \indent
Applying the Riemann-Roch formula for $4P_0 + 2P_1$, we have
\[
 \ell(4P_0 + 2P_1) = \ell(K - 4P_0 - 2P_1) + 1
\]
where $\ell(D) = \dim H^0(X_t, \mathcal{O}(D)) $ and $K$ is the canonical class. 
From the vanishing order of $\omega_i$:
\begin{align*}
\begin{array}{c|cccccc}
& \omega_5 & \omega_3 & \omega_2 & \omega_1 & \omega_6 & \omega_4   
\\ \hline
\text{at} \ P_0 & 0 & 1 & 2 & 3 & 7 & 8
\end{array}, \ \quad
\begin{array}{c|cccccc}
& \omega_5 & \omega_6 & \omega_3 & \omega_4 & \omega_2 & \omega_1    
\\ \hline
\text{at} \ P_1 & 0 & 0 & 1 & 1 & 2 & 3
\end{array},
\end{align*}
we see that there does not exist a holomorphic $1$-form $\omega$ such that 
$\dv(\omega) - 4P_0 - 2P_1$ is positive. Therefore we have $\ell(4P_0 + 2P_1) = 1$ 
and $H^0(X_t, \mathcal{O}(4P_0 + 2 P_1))$ contains only constants functions. 
This implies $\overline{\fA}(4P_0 + 2 P_1) \notin W_{\fA}^5$.
\qed
\subsection{Jacobi inversion} 
We apply Theorem 4 in \cite[Chap. 4, $\S$ 11]{Si71}, for rational functions
\[
 f : X_t \longrightarrow \PP^1,\ (x,y) \mapsto x, \qquad
 g : X_t \longrightarrow \PP^1,\ (x,y) \mapsto 1-x
\]
on $X_t$. Then we have
\begin{align}
\label{J1}
f(Q_1) \times \cdots \times f(Q_6) = \frac{1}{E} \prod_{k=1}^7 \frac{\vartheta
 (\kappa - \fA(Q_1 + \cdots + Q_6) + \int_{\mathfrak{i}_k(\infty,0)} \vec{\xi}, \tau)}
 {\vartheta(\kappa - \fA(Q_1 + \cdots + Q_6) , \tau)}, \\
\label{J2}
g(Q_1) \times \cdots \times g(Q_6) = \frac{1}{E'} \prod_{k=1}^7 \frac{\vartheta
 (\kappa - \fA(Q_1 + \cdots + Q_6) + \int_{\mathfrak{i}_k(\infty,1)} \vec{\xi}, \tau)}
 {\vartheta(\kappa - \fA(Q_1 + \cdots + Q_6), \tau)},
\end{align}
where constants $E$ and $E'$ are independent of $Q_1, \dots, Q_6$, integrals
$\int_{\mathfrak{i}_k(\infty,*)} \vec{\xi} \in \CC^6$ are chosen such that
\[
 \int_{\mathfrak{i}_1(\infty,*)} \vec{\xi} + \cdots + \int_{\mathfrak{i}_7(\infty,*)} \vec{\xi} =0,
\]
 and 
$\fA(Q_1 + \cdots +Q_6) \in \CC^6$ takes the same value in the numerator and the
denominator.
\\ \indent
Substituting $4P_1 + 2P_t$ and $2P_1 + 4 P_t$ for $Q_1 + \cdots + Q_6$ 
in (\ref{J1}), and taking the ratio of resultant equations, we have an expression of 
$t^2$ by theta values:  
\begin{equation} \label{t^2}
\begin{split}
t^2 &= f(P_1)^2 f(P_t)^4 \Big/ f(P_1)^4 f(P_t)^2
\\ 
&= \prod_{k=1}^7 \frac{\vartheta
  (\kappa - \fA(2P_1 + 4P_t) + \int_{\mathfrak{i}_k(\infty,0)} \vec{\xi}, \tau)}
  {\vartheta(\kappa - \fA(2P_1 + 4P_t), \tau)} 
\Big/
\prod_{k=1}^7 \frac{\vartheta
  (\kappa - \fA(4P_1 + 2P_t) + \int_{\mathfrak{i}_k(\infty,0)} \vec{\xi}, \tau)}
  {\vartheta(\kappa - \fA(4P_1 + 2P_t), \tau)} 
\\ 
&= \prod_{k=1}^7 \frac{\vartheta
  (\kappa + a_{4,2} \tau + b_{4,2} + \int_{\mathfrak{i}_k(\infty,0)} \vec{\xi}, \tau)}
  {\vartheta(\kappa + a_{4,2} \tau + b_{4,2}, \tau)} 
\Big/
\prod_{k=1}^7 \frac{\vartheta
  (\kappa + a_{5,2} \tau + b_{5,2} + \int_{\mathfrak{i}_k(\infty,0)} \vec{\xi}, \tau)}
  {\vartheta(\kappa + a_{5,2} \tau + b_{5,2}, \tau)}.
\end{split}
\end{equation}
Similarly, substituting $4P_1 + 2P_t$ and $2P_1 + 4 P_t$ for $Q_1 + \cdots + Q_6$ 
in (\ref{J2}), we have
\begin{equation} \label{1-t^2}
\begin{split}
(1-t)^2 &= g(P_0)^2 g(P_t)^4 \Big/ g(P_1)^4 g(P_t)^2
\\
&=\prod_{k=1}^7 \frac{\vartheta
  (\kappa - \fA(2P_0 + 4P_t) + \int_{\mathfrak{i}_k(\infty,1)} \vec{\xi}, \tau)}
  {\vartheta(\kappa - \fA(2P_0 + 4P_t), \tau)} 
\Big/
  \prod_{k=1}^7 \frac{\vartheta
  (\kappa - \fA(4P_0 + 2P_t) + \int_{\mathfrak{i}_k(\infty,1)} \vec{\xi}, \tau)}
  {\vartheta(\kappa - \fA(4P_0 + 2P_t), \tau)}
\\ 
&=\prod_{k=1}^7 \frac{\vartheta
  (\kappa + a_{2,4} \tau + b_{2,4} + \int_{\mathfrak{i}_k(\infty,1)} \vec{\xi}, \tau)}
  {\vartheta(\kappa + a_{2,4} \tau + b_{2,4}, \tau)} 
\Big/
  \prod_{k=1}^7 \frac{\vartheta
  (\kappa + a_{5,2} \tau + b_{5,2} + \int_{\mathfrak{i}_k(\infty,1)} \vec{\xi}, \tau)}
  {\vartheta(\kappa + a_{5,2} \tau + b_{5,2}, \tau)}.
\end{split}
\end{equation}
\subsection{theta functions with characteristics}
The above expressions are simplified by introducing theta functions with characteristcs 
$a, b \in \QQ^6$:
\begin{align*}
\vartheta_{a,b}(z,\tau) &= \exp[\pi i a\tau {}^ta + 2\pi i a^t(z + b)] 
\vartheta(z + a \tau + b, \tau) \\
&= \sum_{n \in \ZZ^6}\exp[\pi i (n+a) \tau {}^t(n+a) + 2 \pi i (n+a){}^t(z+b)].
\end{align*}
We denote a theta constant $\vartheta_{a,b}(0, \tau)$ by $\vartheta_{a,b}(\tau)$. 
Let $\vartheta_{[m,n]}(z,\tau)$ be $\vartheta_{a,b}(z,\tau)$ with characteristics 
$a = a_{m,n} + a_0,\ b = b_{m,n} + b_0$ in Lemma \ref{Lemma-torsion}. 
With this notation, theta expressions (\ref{t^2}) and (\ref{1-t^2}) are 
\begin{align*}
t^2 = \prod_{k=1}^7 \frac
 {\vartheta_{[2,5]}(\tau) \, \vartheta_{[4,2]}(\int_{\mathfrak{i}_k(\infty,0)} \vec{\xi}, \tau)}
 {\vartheta_{[4,2]}(\tau) \, \vartheta_{[2,5]}(\int_{\mathfrak{i}_k(\infty,0)} \vec{\xi}, \tau)},
\qquad 
(1-t)^2 = \prod_{k=1}^7 \frac
 {\vartheta_{[5,2]}(\tau) \, \vartheta_{[2,4]}(\int_{\mathfrak{i}_k(\infty,1)} \vec{\xi}, \tau)}
 {\vartheta_{[2,4]}(\tau) \, \vartheta_{[5,2]}(\int_{\mathfrak{i}_k(\infty,1)} \vec{\xi}, \tau)}.
\end{align*}
Putting
\[
 \int_{\mathfrak{i}_k(\infty,x)} \vec{\xi} = \begin{cases} a_{1,0} \tau + b_{1,0} \ (x=0) \\ 
a_{0,1} \tau + b_{0,1} \ (x=1) \end{cases} (1 \leq k \leq 6),
 \qquad 
\int_{\mathfrak{i}_7(\infty,x)} \vec{\xi} = \begin{cases} -6(a_{1,0} \tau + b_{1,0}) \ (x=0) \\
-6(a_{1,0} \tau + b_{1,0}) \ (x=1) \end{cases}
\]
and using formulas
\begin{align*}
\vartheta_{a,b}(a' \tau + b', \tau) = 
\exp[-\pi i a' \tau {}^ta' - 2\pi i a' {}^t(b + b')] 
\vartheta_{a+a', b+b'}(0, \tau), \quad a',b' \in \QQ^6, 
\\ 
\theta_{(a + a',b + b')}(z,\tau) = \exp(2 \pi \sqrt{-1}a{}^tb') 
\theta_{(a,b)}(z,\tau), \qquad a', b' \in \ZZ^6,
\end{align*}
we see that 
\begin{align*}
 \prod_{k=1}^7 \frac 
 {\vartheta_{[4,2]}(\int_{\mathfrak{i}_k(\infty,0)} \vec{\xi}, \tau)}
 {\vartheta_{[2,5]}(\int_{\mathfrak{i}_k(\infty,0)} \vec{\xi}, \tau)}
= \z^3 \frac{\vartheta_{[5,2]}(\tau)^7}{\vartheta_{[3,5]}(\tau)^7},
\qquad 
 \prod_{k=1}^7 \frac 
 {\vartheta_{[2,4]}(\int_{\mathfrak{i}_k(\infty,1)} \vec{\xi}, \tau)}
 {\vartheta_{[5,2]}(\int_{\mathfrak{i}_k(\infty,1)} \vec{\xi}, \tau)}
= \frac{\vartheta_{[2,5]}(\tau)^7}{\vartheta_{[5,3]}(\tau)^7}.
\end{align*}
Since $\vartheta_{-a,-b}(-z,\tau) = \vartheta_{a,b}(z,\tau)$, we have
\begin{align*}
 t^2 = \z^3 \frac{ \vartheta_{[5,2]}(\tau)^{14}}{\vartheta_{[4,2]}(\tau)^{14}}, 
\qquad
 (1-t)^2 = \frac{ \vartheta_{[2,5]}(\tau)^{14}}{\vartheta_{[2,4]}(\tau)^{14}},
\end{align*}
namely, there exist constants $\varepsilon_1 = \pm 1$ and $\varepsilon_2 = \pm 1$ 
such that
\begin{align} \label{t and 1-t}
 t = \z^3 \varepsilon_1 \frac{ \vartheta_{[5,2]}(\tau)^7}{\vartheta_{[4,2]}(\tau)^7}, 
\qquad
 1-t = \varepsilon_2 \frac{ \vartheta_{[2,5]}(\tau)^7}{\vartheta_{[2,4]}(\tau)^7}.
\end{align}
\subsection{Theta transformation}
For $g = \begin{bmatrix} A & B \\ C & D \end{bmatrix} \in \mathrm{Sp}_{2g}(\ZZ)$, 
theta constants $\vartheta_{a,b}(\tau)$ satisfy 
the transformation formula  
\begin{align*}
\vartheta_{g(a,b)}(g\tau) = \mu(g) \exp[2 \pi i \lambda_{a,b}(g)]
\mathrm{det}(C\tau + D)^{\frac{1}{2}}\vartheta_{(a,b)}(\tau) 
\end{align*}
where  
\begin{align*}
\lambda_{a,b}(g) &= -\frac{1}{2}(^ta^tDBa - 2^ta^tBCb + {}^tb^tCAb) + 
\frac{1}{2}(^ta^tD - ^tb^tC) \mathrm{diag}(A^tB)
\end{align*}
and $\mu(g)$ is a certain 8-th root of 1 depending only on $g$. 
Therefore, as coordinates of $\PP^2(\CC)$, we have
\begin{align} \label{proj}
 [\vartheta_{g[2,4]}:\vartheta_{g[2,5]}:\vartheta_{g[3,5]}](g \cdot \tau) 
=[\ee[\lambda_{[2,4]}(g)] \vartheta_{[2,4]} : \ee[\lambda_{[2,5]}(g)] \vartheta_{[2,5]} : 
\ee[\lambda_{[3,5]}(g)] \vartheta_{[3,5]}](\tau) 
\end{align} 
where $\ee[-] = \exp[2 \pi i -]$. 
\\ \indent
By explicit form of $\sigma_0 = \phi(h_0)$ and $\sigma_1 = \phi(h_1)$ in Appendix, 
we see that
\begin{align*}
 \lambda_{[2,4]}(\sigma_0) &= 53/56, & 
 \lambda_{[2,5]}(\sigma_0) &= 53/56, & 
 \lambda_{[3,5]}(\sigma_0) &= 7/8, 
\\
 \lambda_{[2,4]}(\sigma_1) &= 25/56, & 
 \lambda_{[2,5]}(\sigma_1) &= 19/392, & 
 \lambda_{[3,5]}(\sigma_1) &= 79/392,
\end{align*}
and 
\begin{align*}
\vartheta_{\sigma_0[2,4]} &= \ee[5/14] \vartheta_{[5,2]}, &
\vartheta_{\sigma_0[2,5]} &= - \vartheta_{[5,3]}, &
\vartheta_{\sigma_0[3,5]} &= \ee[13/14] \vartheta_{[4,2]}
\\
\vartheta_{\sigma_1[2,4]} &= -\vartheta_{[5,3]}, &
\vartheta_{\sigma_1[2,5]} &= \ee[9/14] \vartheta_{[4,2]}, &
\vartheta_{\sigma_1[3,5]} &= \ee[4/7] \vartheta_{[5,2]}. 
\end{align*}
Applying these for (\ref{proj}), we obtain 
\begin{equation} \label{modular-action}
\begin{split}
 [\vartheta_{[2,4]} : \vartheta_{[2,5]} : \vartheta_{[3,5]}](\sigma_0 \cdot \tau)  
&=[-\vartheta_{[2,5]} : \ee[9/14] \vartheta_{[2,4]} : \vartheta_{[3,5]}](\tau),
\\
 [\vartheta_{[2,4]} : \vartheta_{[2,5]} : \vartheta_{[3,5]}](\sigma_1 \cdot \tau)  
&=[\vartheta_{[2,4]} : \ee[67/98] \vartheta_{[3,5]} : \ee[45/98] \vartheta_{[2,5]}](\tau).
\end{split}
\end{equation}
\subsection{Theorem} \label{main}
(1) The inverse of the Schwarz map
\begin{align*}  
 \mathfrak{s} : \CC-\{0,1\} \longrightarrow \DD_H^+, \qquad 
 t \mapsto u = [u_1(t):u_3(t)]
\end{align*}
is given by $\Gamma(1-\z)$-invariant function $\mathfrak{t}(u) = 
\z^5 \frac{ \vartheta_{[2,5]}(\Phi(u))^7}{\vartheta_{[3,5]}(\Phi(u))^7}$, where 
$\Phi : \DD_H^+ \rightarrow \HH_6^M$ is the modular embedding given in Appendix. 
In other word, $\Phi(u) \in \HH_6^M$ is the period matrix of an algebraic curve 
\[
 y^7 = x(x-1)(x-\mathfrak{t}(u)).
\]
(2) The analytic map 
\begin{align*}
 Th \ : \ \DD_H^+ \longrightarrow \PP^2(\CC), \qquad u \mapsto 
  [\ee[5/49] \vartheta_{[2,4]} \vartheta_{[2,5]} : \vartheta_{[2,5]} \vartheta_{[3,5]} : 
  -\vartheta_{[2,4]} \vartheta_{[3,5]}](\Phi(u))
\end{align*}
induces an isomorphism $\DD_H^+ / \Gamma((1-\z)^2)$ and the Fermat septic curve
\begin{align*}
 \mathcal{F}_7 \ : \ X^7 + Y^7 + Z^7 = 0, \qquad [X:Y:Z] \in \PP^2(\CC).
\end{align*}
\noindent
{\it Proof}. 
From (\ref{t and 1-t}), we have 
\begin{align*}
1 = \varepsilon_1 \z^5 \frac{ \vartheta_{[2,5]}(\tau)^7}{\vartheta_{[3,5]}(\tau)^7} 
+ \varepsilon_2 \frac{ \vartheta_{[2,5]}(\tau)^7}{\vartheta_{[2,4]}(\tau)^7}.
\end{align*}
Since this equation must be invariant under actions of $\sigma_0 = \phi(h_0)$ and 
$\sigma_1 = \phi(h_1)$ in (\ref{modular-action}) (otherwise, the image of $Th$ 
is not irreducible), we see that $\varepsilon_1 = \varepsilon_2 = 1$ and 
\[
 t = \z^5 \frac{ \vartheta_{[2,5]}(\Phi(u))^7}{\vartheta_{[3,5]}(\Phi(u))^7}. 
\] 
Let us recall that $\Gamma(1-\z)$ is projectively generated by $h_0^2$ and $h_1^2$, 
and $\Gamma((1-\z)^2)$ is projectively isomorphic to the commutator subgroup of
$\Gamma(1-\z)$. From (\ref{modular-action}), we see that
\begin{align*}
[\vartheta_{[2,4]} : \vartheta_{[2,5]} : \vartheta_{[3,5]}](\sigma_0^2 \cdot \tau)  
=[\z \vartheta_{[2,4]} : \z \vartheta_{[2,5]} : \vartheta_{[3,5]}](\tau),
\\
 [\vartheta_{[2,4]} : \vartheta_{[2,5]} : \vartheta_{[3,5]}](\sigma_1^2 \cdot \tau)  
=[\vartheta_{[2,4]} : \z \vartheta_{[2,5]} : \z \vartheta_{[3,5]}](\tau).
\end{align*}
Therefore the commutator subgroup of $\Gamma(1-\z)$ acts trivialy on 
\[
 [\vartheta_{[2,4]}(\Phi(u)) : \vartheta_{[2,5]}(\Phi(u)) : \vartheta_{[3,5]}(\Phi(u))] 
\in \PP^2,
\] 
and the map $Th$ gives a $(\ZZ/7\ZZ)^2$-equivariant isomorphism of 
$\DD_H^+ / \Gamma((1-\z)^2)$ and the Fermat septic curve.
\qed
\subsection{Klein quartic}
It is known that the Klein quartic curve  
\[
 \mathcal{K}_4 \ : \ X^3 Y + Y^3 Z + Z^3 X = 0, \qquad [X:Y:Z] \in \PP^2(\CC).
\]
is the quotient of $\mathcal{F}_7$ by an automorphism
\[
 \alpha \ : \ \mathcal{F}_7 \longrightarrow \mathcal{F}_7, \qquad 
 [X:Y:Z] \mapsto [\z X: \z^3 Y:Z]
\]
which is induced by $g_0 g_1^3 \in \Gamma(1-\z)$ via the map $Th$. 
The quotient map is given by 
\[
 \mathcal{F}_7 \longrightarrow \mathcal{K}_4, \qquad
 [X:Y:Z] \mapsto [XY^3 : YZ^3 : ZX^3].
\]
The Klein quartic $\mathcal{K}_4$ is isomorphic to the elliptic modular curve 
$\mathcal{X}(7)$ of level $7$, and also to a Shimura curve parametrizing a family of 
QM Abelian $6$-folds (see \cite{E99}). The following Corollary gives a new moduli  
interpretation of $\mathcal{K}_4$.   
\subsection{Corollary} \label{Klein1}
The Klein quartic curve $\mathcal{K}_4$ is isomorphic to $\DD_H^+ / \Gamma_{Klein}$ 
where 
\[
 \Gamma_{Klein} = \{ \begin{bmatrix} a & b \\ c & d \end{bmatrix} 
 \in \Gamma(1-\z) \ | \ a \equiv 1 \mod (1-\z)^2 \}. 
\]
{\it Proof}.
Let us recall the homomorphism 
\begin{align*}
  \nu : \Gamma(1-\z) \longrightarrow \mathrm{M}_2(\FF_7), 
\qquad
 \nu(g) = \frac{1}{1-\z} (g - 1) \mod 1-\z
\end{align*}
in the proof of Proposition \ref{Prop-group}. The kernel of $\nu$ is $\Gamma((1-\z)^2)$
 and the image is generated by
\begin{align*}
  \nu(g_0) = \begin{bmatrix} -1 & 0 \\ 0 & 0 \end{bmatrix}, \quad
    \nu(g_1) = \begin{bmatrix} 5 & 1 \\ 5 & 1 \end{bmatrix}.
\end{align*}
Since we have $\nu(g_0^a g_1^b) = \begin{bmatrix} -a + 5b & b \\ 5b & b \end{bmatrix}$,  
the group $\Gamma_{Klein}$ is generated by $\Gamma((1-\z)^2)$ and $g_0 g_1^3$. Namely 
we have $\DD_H^+ / \Gamma_{Klein} = \mathcal{F}_7 / \left< \alpha \right>$.  
\qed
\subsection{PEL-family} 
Let $(A, E, \rho, \lambda)$ be a $4$-tuple
\begin{itemize}
\item[(1)] $A$ is a $6$-dimensional complex Abelian variety $V/\Lambda$, where 
$V$ is ismorphic to the tangent space $T_0A$ and $\Lambda$ is isomorphic to $H_1(A,\ZZ)$. 
\item[(2)] $E \ : \ \Lambda \times \Lambda \rightarrow \ZZ$ is 
a principal polarization.
\item[(3)] $\rho$ is an automorphism of order $7$ preserving $E$, and the induced action 
on $T_0A$ has eigenvalues $\z, \z, \z^2, \z^2, \z^3, \z^4$.
\item[(4)] $\lambda : \Lambda \rightarrow \ZZ[\z]^2$ is an isomorphism such that 
\[
\lambda(\rho(x)) = \z^4 \lambda(x), \qquad
E(x,y) = \frac{1}{7} Tr_{\QQ(\z)/\QQ}
((\z^3-\z^4) ^t \overline{\lambda(x)} H^{-1} \lambda(y))
\] 
(see Remark \ref{Rem Riemann form}). Note that $\lambda$ induces an isomorphism of 
the torsion subgroup $A_{tor}$ and $(\QQ(\z) / \ZZ[\z])^2$.
\end{itemize}
An isomorphism $f : (A, E, \rho, \lambda) \rightarrow (A', E', \rho', \lambda')$ is 
defined as an isomorphism of Abelian varieties $f:A \rightarrow A'$ such that 
$f^*E' = E$, $f \circ \rho = \rho' \circ f$ and $\lambda = \lambda' \circ f$. 
Then we see that
\subsection{Corollary} \label{Klein2}
We have isomrophisms
\begin{align*}
 \DD_H^+ / \Gamma(\mathfrak{m}) &\cong
 \begin{Bmatrix}
   \text{Set of} \ (A, E, \rho, \lambda) \ \text{modulo isomorphisms} 
   \ $f$ \ \text{such that} \\
  \lambda^{-1} \equiv (\lambda' \circ f)^{-1} \ \ \text{on} \ \
(\mathfrak{m}^{-1} \ZZ[\z] / \ZZ[\z])^2 
 \end{Bmatrix},
\\
\DD_H^+ / \Gamma_{Klein} &\cong
 \begin{Bmatrix}
   \text{Set of} \ (A, E, \rho, \lambda) \ \text{modulo isomorphisms} 
   \ $f$ \ \text{such that} \hspace{2cm} \\
  \text{(i)} \ \lambda^{-1}(\frac{1}{(1-\z)^2},0) = 
(\lambda' \circ f)^{-1}(\frac{1}{(1-\z)^2}, \frac{b}{1-\z}) \ \text{for} \ 
{}^{\exists} b \in \ZZ[\z] \\
  \text{(ii)} \ \lambda^{-1}(0, \frac{1}{1-\z}) = 
(\lambda' \circ f)^{-1}(0, \frac{1}{1-\z}) \hspace{3cm}
 \end{Bmatrix}.
\end{align*} 
\section{K3 surface}
\subsection{}
In this final section, we construct K3 surfaces with a non-symplectic automorphism 
of order $7$ attached to $X_t$, according to Garbagnati and Penegini (\cite{GP}). 
For generalities on K3 surfaces and elliptic surfaces, see \cite{SS10} and references 
in it. Let us consider two curves
\begin{align*}
X_t \ : \ y_1^7 = x_1(x_1 -1)(x_1 - t), \qquad
X_{\infty} \ : \ y_2^7 = x_2^2 - 1 
\end{align*}
and an affine algebraic surface
\[
S_t \ : y^2 = x(x-z)(x-tz)+z^{10}.
\]
$X_t$ is a hyperelliptic curve of genus $3$. The surface $S_t$ is birational to the 
quotient of $X_t \times X_{\infty}$ by an automorphism
\begin{align*}
\rho \times \rho \ : \ X_t \times X_{\infty} \longrightarrow X_t \times X_{\infty},
\qquad
 (x_1, y_1) \times (x_2, y_2) \mapsto (x_1, \z y_1) \times (x_2, \z y_2), 
\end{align*}
and the ratioinal quotient map $X_t \times X_{\infty} \dashrightarrow S_t$ is given by
\begin{align*}
 z = y_1 / y_2, \qquad  y = z^5 x_2, \qquad x = z x_1.
\end{align*}
The minimal smooth compact model of $S_t$ (denoted by the same symbol $S_t$) is a K3 
surface with an elliptic fibration
\[
 \pi \ : \ S_t \longrightarrow \PP^1, \quad (x,y,z) \mapsto z.
\]
To see this, let us consider a minimal Weierstrass form 
\begin{align*}
S'_t \ : \ y^2 = &x^3 + G_2(z) x + G_3(z) \\
  &G_2(z) = -\frac{1}{3}(t^2 - t +1) z^2, \qquad
  G_3(z) = z^{10} - \frac{1}{27}(2t-1)(t+1)(t-2) z^3 
\end{align*}
and the discriminant
\begin{align*}
\Delta(z) = 4G_2(z)^3 + 27 G_3(z)^2
 = z^6 \{ 27z^{14} -2 (2t-1)(t+1)(t-2) z^7 - t^2(t-1)^2 \}.
\end{align*}
From this, we see that $S_t$ is a K3 surface, and it has a singular fiber of type 
$\mathrm{I}_0^*$ at $z=0$, of type $\mathrm{IV}$ at $z = \infty$ and fourteen fibers of 
type $\mathrm{I}_1$ on $\PP^1 - \{ 0, \infty \}$. 
Note that 
\begin{align*}
 \frac{dx_1}{y_1^3} \otimes \frac{y_2^2 d y_2}{x_2} \in 
 H^0(X_t, \Omega^1) \otimes H^0(X_{\infty}, \Omega^1)
\end{align*}
is the unique $(\rho \times \rho)$ - invariant element up to constants,  
and descents to a holomorphic $2$-form on $S_t$ (see \cite{GP}, Section 3). Therefore 
the period map for a family of K3 surface $S_t$ is given by the Schwarz map $\mathfrak{s}$.  
Note also that an automorphism $\rho \times \mathrm{id}$ of $X_t \times X_{\infty}$ descentds 
to $S_t$:
\[
 \rho \times \mathrm{id} \ : \ S_t \longrightarrow S_t, \quad (x,y,z) \mapsto 
(\z x,\z^5 y,\z z).
\]
Since $S_t / \left< \rho \times \mathrm{id} \right>$ is birational to a rational surface
$X_t / \left< \rho \right> \times X_{\infty} / \left< \rho \right>$, the automorphism 
$\rho \times \mathrm{id}$ is non-simplectic. Hence the transcendental lattice 
$\mathrm{T}_{S_t}$ is a free $\ZZ[\rho \times \mathrm{id}]$-module (\cite{Ni79}). 
Since our family has positive dimensional moduli, we have 
$\mathrm{rank} \: \mathrm{T}_{S_t} \geq 12$ and $\mathrm{rank} \: \mathrm{NS}(S_t) \leq 10$ 
for a general $t \in \CC-\{0,1\}$, where $\mathrm{NS}(S_t)$ is the N\'eron-Severi lattice. 
\subsection{} 
Let us compute the N\'eron-Severi lattice and the Mordell-Weil group $\mathrm{MW}(S_t)$. 
Let $o$ be the zero section of $\pi : S_t \rightarrow \PP^1$. We have three sections 
\[
 s_a \ : \ \PP^1 \longrightarrow S_t, \quad z \mapsto (x,y,z)=(az,z^5,z), 
\qquad a =0,\: 1,\: t
\]
such that $s_0 + s_1 + s_t = o$ in $\mathrm{MW}(S_t)$.
Let $2 \ell_0 + \ell_1 + \ell_2 + \ell_3 + \ell_4$ be the irreducible decomposition 
of $\pi^{-1}(0)$, and $\ell_1' + \ell_2' + \ell_3'$ be that of $\pi^{-1}(\infty)$. 
For a suitable choice of indeces, intersection numbers of these curves are given by
 the following graph; self intersection number of each curve is $-2$, two curves are 
connected by an edge if they intersect and intersection numbers are $1$ except 
$s_a \cdot s_b = 2$.
\begin{center}
\begin{picture}(260,130)(0,0) 
\put(0,60){$\ell_0$} \put(4,63){\circle{13}} 
\put(60,0){$\ell_1$} \put(60,40){$\ell_2$} \put(60,80){$\ell_3$} \put(60,120){$\ell_4$}
\multiput(64,3)(0,40){4}{\circle{13}}
\put(120,0){$o$} \put(120,40){$s_0$} \put(120,80){$s_1$} \put(120,120){$s_t$}
\multiput(124,3)(0,40){4}{\circle{13}}
\put(180,80){$\ell'_2$} \put(184,83){\circle{14}}
\put(180,0){$\ell'_1$} \put(184,3){\circle{14}}
\put(240,40){$\ell'_3$} \put(244,43){\circle{14}}
\put(7,69){\line(1,1){51}}
\put(10,66){\line(5,2){48}}
\put(10,60){\line(5,-2){48}}
\put(7,57){\line(1,-1){51}}
\multiput(71,3)(0,40){4}{\line(1,0){46}}
\put(131,3){\line(1,0){46}} \put(131,83){\line(1,0){46}}
\put(130,46){\line(5,3){51}} \put(130,120){\line(5,-3){51}}
\put(184,10){\line(0,1){66}}
\put(191,5){\line(5,3){52}} \put(191,81){\line(5,-3){52}}
\qbezier(118,46)(95,80)(118,120)
\thicklines
\put(124,89){\line(0,1){27}} \put(124,49){\line(0,1){27}}
\put(125,63){$2$} \put(125,97){$2$} \put(100,95){$2$} 
\end{picture}
\end{center}
Let $N \subset \mathrm{NS}(S_t)$ be the lattice generated by 
$o, \ s_0, \ s_1, \ s_t, \ \ell_0, \ \ell_1, \ \ell_2, \ \ell_3, \ \ell_4, \ \ell_1'$. 
The rank of $N$ is $10$ and the discriminant is $-49$. Hence the Picard number of $S_t$ 
is generically $10$ and the rank of $\mathrm{MW}(S_t)$ is $2$ by the Shioda-Tate formula 
(\cite{SS10}, Corollary 6.13). Since the fixed locus $S_t^{\rho \times \mathrm{id}}$ is contained
 in $\pi^{-1}(0) \cup \pi^{-1}(\infty)$ and no elliptic curve contained in 
$S_t^{\rho \times \mathrm{id}}$, we see that 
$\mathrm{NS}(S_t) = \mathrm{U}(7) \oplus \mathrm{E}_8$ by the classification theorem of
 Artebani, Sarti and Taki (\cite{AST11}, $\S$ 6). Therefore we have $\mathrm{NS}(S_t) = N$. 
Let $L$ be the lattice generated by the zero section and vertical divisors. It is known that 
$\mathrm{MW}(S_t) \cong \mathrm{NS}(S_t)/L$ (\cite{SS10}, Theorem 6.3). Now it is obvious 
that $\mathrm{MW}(S_t) = \ZZ s_0 \oplus \ZZ s_1 \cong \ZZ^2$. 
\appendix
\section{}
\subsection{Symplectic representatioin}
\[
M =
\footnotesize
\left[ \begin{array}{ccc|ccc|ccc|ccc}
 0 & 0 & 0 &  &  &  & -1 & 1 & 0 &  &  &  \\
 0 & 0 & 0 &  & \bf{O} &  & 0 & -1 & 1 &  & \bf{O} &  \\
 -1 & -1 & -1 &  &  &  & 0 & 0 & -1 &  &  &  \\
\hline
  &  &  & -1 & 0 & 1 &  &  &  & 0 & 1 & 0 \\
  & \bf{O} &  & 0 & 0 & -2 &  & \bf{O} &  & 1 & -1 & 1 \\
  &  &  & 0 & -1 & 1 &  &  &  & 0 & 1 & -1 \\
\hline
 1 & 0 & 0 &  &  &  &  &  &  &  &  &  \\
 1 & 1 & 0 &  & \bf{O} &  &  & \bf{O} &  &  & \bf{O} &  \\
 1 & 1 & 1 &  &  &  &  &  &  &  &  &  \\
\hline
  &  &  & 1 & -1 & -2 &  &  &  & 0 & -1 & 0 \\
  & \bf{O} &  & -1 & 1 & 1 &  & \bf{O} &  & -1 & 0 & 0 \\
  &  &  & -1 & 0 & 3 &  &  &  & -1 & 1 & -1 
\end{array} \right]
\]
\normalsize
\[
\phi(h_0) =
\footnotesize
\left[ \begin{array}{ccc|ccc|ccc|ccc}
 0 & 0 & 0 &  &  &  & 1 & -1 & 0 & & & \\
 0 & 0 & 0 &  & \bf{O} &  & 0 & 1 & -1 & & \bf{O}& \\
 1 & 1 & 1 &  &  &  & 0 & 0 & 1 & & & \\
\hline
 & & & 1 & 0 & 0 & & & & & & \\
 & \bf{O} & & 0 & 1 & 0 & & \bf{O} & & & \bf{O} & \\
 & & & 0 & 0 & 1 & & & & & & \\
\hline
 -1 & 0 & 0 & & & & & & & & & \\
 -1 & -1 & 0 & & \bf{O} & & & \bf{O} & & & \bf{O} & \\
 -1 & -1 & -1 & & & & & & & & & \\
\hline
 & & & & & & & & & 1 & 0 & 0 \\
 & \bf{O} & & & \bf{O} & & & \bf{O} & & 0 & 1 & 0 \\
 & & & & & & & & & 0 & 0 & 1 
\end{array} \right] 
\]
\normalsize
\[
\phi(h_1) =
\footnotesize
\left[ \begin{array}{ccc|ccc|ccc|ccc}
 1 & 1 & 1 & 0 & 0 & 0 & 1 & 0 & 0 & 1 & 1 & 0 \\
 0 & 1 & 1 & 0 & 1 & 0 & 0 & 1 & 0 & 1 & 0 & 1 \\
 0 & 0 & 1 & 1 & 0 & 0 & 0 & 0 & 1 & 0 & 0 & 0 \\
\hline
 0 & 1 & 2 & 1 & 1 & 0 & 1 & 1 & 0 & 2 & 1 & 1 \\
 1 & 2 & 2 & 0 & 1 & 2 & 0 & 1 & 1 & 1 & 2 & 0 \\
 0 & 0 & 1 & 0 & 1 & 0 & 1 & 0 & 0 & 1 & 0 & 1 \\
\hline
 -1 & -1 & -1 & 0 & 0 & 1 & 0 & 0 & 0 & -1 & 0 & -1 \\
 -1 & -2 & -2 & 0 & 0 & 0 & -1 & 0 & 0 & -2 & -1 & -1 \\
 -1 & -2 & -3 & 0 & -1 & -1 & -1 & -1 & 0 & -2 & -2 & -1 \\
\hline
 1 & 1 & 1 & -1 & 0 & 1 & 0 & 0 & 0 & 1 & 1 & 0 \\
 0 & 0 & -1 & 0 & -1 & -1 & 0 & 0 & -1 & 0 & 0 & 0 \\
 -1 & -1 & -1 & 0 & 0 & -2 & 1 & -1 & -1 & 0 & -1 & 1 \\
\end{array} \right] 
\]
\normalsize
\subsection{Period matrix} Let $z$ be $\exp[2 \pi i/7]$ and put
\footnotesize
\begin{align*}
a = 1 + z + z^2 + z^4 = \frac{1 + \sqrt{-7}}{2}, \qquad
b_1 = z - 2z^2 - 2z^4, \qquad  b_2 = -(2 z^3 + 1 - z^6 + 2 z^5).
\end{align*}
\normalsize
The modular embedding $\Phi : \DD_H^+ \rightarrow \HH_6^M$ in Proposition 
\ref{modular-embedding} is given by
\begin{align*}
  \Phi(u) = \frac{1}{\Delta}
\left( \begin{bmatrix} A_{11} & O \\ O & D_{11} \end{bmatrix} u_1^2 
+ \begin{bmatrix} O & B_{12} \\ {}^tB_{12} & O \end{bmatrix} u_1 u_2 
+ \begin{bmatrix} A_{22} & O \\ O & D_{22} \end{bmatrix} u_2^2 \right)
\end{align*}
where 
\[
\Delta = (z^2+z+1)(2 z^2-z+2) u_1^2 + 3 (z+1) u_2^2
\]
and
\footnotesize
\begin{align*}
A_{11} &=  
(z^2 + z + 1)(2 z^2 - z + 2) 
\left[ \begin{array}{ccc}
 a & 0 & -1 \\
 0 & a-1 & -a \\
 -1 & -a & 1 
\end{array} \right],
\quad 
D_{11} = (z^2 + 1) 
\left[ \begin{array}{ccc}
 2 z^6+z^5-z^3-1 & 2 z^6-z^3 & -z^3 \\
 2 z^6-z^3 & z^2-z^3 & z^6 - a \\
 -z^3 & z^6 - a & a 
\end{array} \right]
\\
B_{12} &= (z^3 - z^5) 
\left[ \begin{array}{ccc}
 -b_1 & -b_2 & -1 \\
 z^5 b_1 & z^5 b_2 & z^5 \\
 (1 + z^6) b_1 & (1 + z^6) b_2 & (1 + z^6) \\
\end{array} \right]
\\
A_{22} &= -3 
\left[ \begin{array}{ccc}
 z(z^5 - z^2 - 1) & z-1 & z^3 + 1 \\
 z-1 & z^5-z^2+1 & z^2 \\
 z^3 + 1 & z^2 & -z^2 (z+1) 
\end{array} \right],
\quad
D_{22} = (z+1) 
\left[ \begin{array}{ccc}
 3 a -2  & a-1 & a \\
 a-1 & 2 a - 1 & -2 \\
 a & -2 & a + 1 
\end{array} \right]. 
\end{align*}


\begin{thebibliography}{35}
\bibitem[AST11]{AST11}
 M. Artebani, A. Sarti and S. Taki, 
        {\it K3 surfaces with non-symplectic automorphisms of prime order}, 
        Math. Z., {\bf 268} (2011), 507-533.

\bibitem[BR87]{BR87}
 M. Bershadsky and A. Radul, 
        {\it Conformal field theories with additional $Z_N$ symmetry},
        Int. J. Mod. Phys., {\bf A2-1} (1987), 165-178.

\bibitem[BR88]{BR88}
 M. Bershadsky and A. Radul, 
        {\it Fermionic fields on $Z_N$ curves}, 
        Commun. Math. Phys. {\bf 116} (1988), no. 4, 689-700.

\bibitem[B07]{B07}
 F. Beukers, {\it Gauss' hypergeometric function}, in {\it Arithmetic and Geometry Around 
        Hypergeometric Functions}, Progress in Math., vol. {\bf 260}, 23-42. 
        Birkh\"auser (2007).

\bibitem[BL92]{BL92}
 Ch. Birkenhake and H. Lange, {\it Complex Abelian Varieties}, Springer (1992).

\bibitem[CIW94]{CIW94}
 P. B. Cohen, C. Itzykson and J. Wolfart, {\it Fuchsian triangle groups and Grothendieck 
        dessins. Variations on a theme of Belyi}, Comm. Math. Phys., {\bf 163} (1994), 
        no. 3, 605-627.  

\bibitem[DGMS13]{DGMS13}
 C. Doran, T. Gannon, H. Movasati and K. Shokri, {\it Automorphic forms for triangle groups},
        Commun. Number Theory Phys. {\bf 7} (2013), no. 4, 689-737. 

\bibitem[E99]{E99}
 N. D. Elkies, {\it The Klein quartic in number theory}, in {\it The Eightfold Way}, 
       Math. Sci. Res. Inst. Publ., {\bf 35}, 51-102, Cambridge Univ. Press (1999).

\bibitem[EG06]{EG06}
 V.Z. Enolski AND T.Grava, {\it Thomae type formulae for singular $Z_N$ curves}, 
       Lett. Math. Phys. {\bf 76} (2006), no. 2-3, 187-214. 

\bibitem[GP]{GP}
A. Garbagnati and M. Penegini, {\it K3 surfaces with a non-symplectic automorphism and 
        product-quotient surfaces with cyclic group}, arXiv1303.1653[math.AG].

\bibitem[vG92]{vG92} 
 B. van Geemen, {\it Projective models of Picard Modular Varieties}, 
        in {\it Classification of Irregular Varieties}, Springer LNM {\bf 1515} (1992),
        68-99.

\bibitem[HY99]{HY99}
 M. Hanamura and M. Yoshida, {\it Hodge structure on twisted cohomologies}, 
        Nagoya Math. J., Vol. {\bf 154} (1999), 123-139.

\bibitem[H05]{H05} 
 M. Harmer, {\it Note on the Schwarz triangle functions}, Bull. Austral. Math. Soc. 
        {\bf 72} (2005), no. 3, 385-389. 

\bibitem[dJN91]{dJN91}
 J. de Jong and R. Noot, {\it Jacobians with complex multiplication}, in
        {\it Arithmetic Algebraic Geometry}, Progress in Math., vol. {\bf 89}, pp. 177-192,
        Birkh\"auser (1991). 

\bibitem[K03]{K03}
 K. Koike, {\it On the family of pentagonal curves of genus 6 and associated modular 
        forms on the ball}, J. Math. Soc. Japan, {\bf 55} (2003), 165 - 196.

\bibitem[KS07]{KS07}
 K. Koike and H. Shiga, {\it Isogeny formulas for the Picard modular form and a three 
        terms arithmetic geometric mean}, J. Number Theory, {\bf 124} (2007), 123-141.

\bibitem[KW04]{KW04}
 K. Koike and A. Weng, {\it Construction of CM Picard curves}, Math. Comp., {\bf 74}
        (2004), 499-518.

\bibitem[Mi75]{Mi75} 
 J. Milnor, {\it On the 3-dimensional Brieskorn manifolds $M(p,q,r)$}, in 
        {\it Knots, groups, and 3-manifolds}, Ann. of Math. Studies, No. {\bf 84}, 
        Princeton Univ. Press (1975), pp. 175-225. 

\bibitem[MO13]{MO13} 
 B. Moonen and F. Oort, {\it The Torelli locus and special subvarieties}, 
        in {\it Handbook of Moduli vol. II}, Adv. Lect. Math., {\bf 25}, 
        Int. Press (2013), 549-594.   

\bibitem[M83]{M83}
 D. Mumford, {\it Tata Lectures on Theta I}, Birkh\"auser (1983).

\bibitem[Na97]{Na97}
 A. Nakayashiki, {\it On the Thomae formula for $Z_N$ curves}, 
        Publ. Res. Inst. Math. Sci. {\bf 33} (1997), 987-1015.

\bibitem[Ni79]{Ni79}
 V. V. Nikulin, {\it Finite automorphism groups of Kahler surfaces of type K3}, 
        Proc. Moscow Math. Soc. {\bf 38} (1979), 75-137.

\bibitem[P1883]{P1883}
 E. Picard, {\it Sur les fonctions de deux vanables ind\'ependantes analogues aux 
        fonctions modulaires}, Acta math., {\bf 2} (1883), 114-135.

\bibitem[R09]{R09}
 C. Rohde, {\it Cyclic Coverings, Calabi-Yau Manifolds and Complex
        Multiplication}, Springer LNM {\bf 1975} (2009).

\bibitem[SS10]{SS10}
 M. Sch\"{u}tt and T. Shioda, {\it Elliptic surfaces}, Algebraic Geometry in 
        East Asia - Seoul 2008, Adv. Stud. Pure Math. {\bf 60} (2010), 51-160.

\bibitem[Si71]{Si71}
 C. L. Siegel, {\it Topics in Complex Function Theory, Vol. II}, Wiley (1971).

\bibitem[Sh79,81]{Sh79-81}
 H. Shiga, {\it On attempt to the K3 modular function I-II}, 
        Ann. Scuola Norm. Sup. Pisa Cl. Sci. (4) {\bf 6} (1979), no. 4, 609-635, 
        {\bf 8} (1981), no. 1, 157-182.

\bibitem[Sh88]{Sh88}
 H. Shiga, {\it On the representation of the Picard modular function by
	$\theta$ constants I-II}, Pub. R.I.M.S. Kyoto Univ. {\bf 24} (1988),
	311-360.

\bibitem[Sm64]{Sm64}
 G. Shimura, {\it On purly transcendental fields of automorphic functions of
	several complex variables}, Osaka J. Math. {\bf 1} (1964), 1-14.

\bibitem[T77]{T77}
 K. Takeuchi, {\it Arithmetic triangle group}, J. Math. Soc. Japan,
	Vol. {\bf 29}, No. 1 (1977), 91-106.

\bibitem[YY84]{YY84}
 T. Yamazaki and M. Yoshida, {\it On Hirzebruch's Examples of Surfaces with
	$c_1^2 = 3c_2$}, Math. Ann. {\bf 266} (1984), 421-431.

\bibitem[Y97]{Y97}
 M. Yoshida, {\it Hypergeometric functions, my love}, Aspects of Mathematics,
        E32. Friedr. Vieweg \& Sohn (1997).

\bibitem[W81]{W81}
 J. Wolfart, {\it Graduierte algebren automorpher formen zu dreiecksgruppen},
        Analysis {\bf 1}, no. 3 (1981), 177-190.
\end{thebibliography}
\end{document}